\documentclass[12pt, a4paper]{amsart}

\usepackage{amsfonts,amsmath,amssymb,amscd}  \usepackage[all]{xy}

\def\Z{{\Bbb Z}}

\def\F{{\Bbb F}}
\def\G{{\Bbb G}}
 
\def\Im{{\mathrm{Im}}}

\def\Ker{{\mathrm{Ker}}}

\def\Im{\mathrm{Im}}

\def\Spec{{\mathrm{Spec}}}
\def\gr{\mathrm{gr}}

\newtheorem{theorem}{Theorem}[section]

\newtheorem{prop}[theorem]{Proposition}
\newtheorem{corollary}[theorem]{Corollary}

\theoremstyle{definition}
\newtheorem{definition}[theorem]{Definition}

\newtheorem{notation}[theorem]{Notation}

\usepackage{color}

\newcommand\pf{\begin{proof}}
\newcommand\epf{\end{proof}}

\usepackage{mathrsfs}

\title[]
{On the Grothendieck resolution for a certain finite flat commutative group scheme of order $p^{n}$ over an $\Bbb{F}_{p}$-algebra\\
}

\author[Yuji~Tsuno]{Yuji~Tsuno}
\address{Yuji Tsuno: National Institute of Technology, Wakayama College,
77 Noshima, Nada-cho, Gobo, Wakayama, Japan 644-0023}
\email{tsuno@wakayama-nct.ac.jp}

\begin{document} 

\begin{abstract}
For any commutative finite flat group scheme,  \\    Grothendieck constructed an embedding of it into some smooth group scheme. This embedding is named the Grothendieck resolution. Let $p$ be a prime number and $n$ a positive integer. In connection with the normal basis problems in the framework of group schemes proposed by Suwa and the author, we consider the Grothendieck resolution for a certain finite flat commutative group scheme of order $p^{n}$ over an $\Bbb{F}_{p}$-algebra.
\end{abstract}

\maketitle

\noindent

\medskip

\section{Introduction}\label{sec:intro}
\renewcommand{\thefootnote}{} 
\footnote[0]{\sc Mathematics Subject Classification (2020). Primary 13B05; Secondary 14L15, 12G05\\
 \ \ \  Keywords and phrases. Finite flat group scheme, normal basis problem, Hopf-Galois extension, cleft extension
}
\renewcommand{\thefootnote}{\arabic{footnote}} 
Let $k$ be a field and $\varGamma$ a finite group. Then, consider the problem, ``Does there exist a Galois extension of $k$ with the Galois group $\varGamma$? If it exists, construct it concretely.''
This is named the Galois inverse problem and is important in Galois theory. The Kummer and Artin-Schreier-Witt theories provide concise descriptions of Galois extensions with a cyclic group. They are important items in Galois theory, with Kummer theory particularly noted for its elementary proof using Lagrange resolvents. 
 Considering the algebraic group $U(\varGamma)$ representing the unit group of the group algebra $k[\varGamma]$, the embedding $\varGamma \rightarrow U(\varGamma)$, and the normal basis theorem, Serre [8, Ch.VI, 8] formulated this method of proof as follows:
Any Galois extension $K/k$ with the Galois group $\varGamma$ is obtained by a Cartesian diagram 
\[\begin{CD}
 \Spec \ K  @>>> U(\varGamma) \\
 @VVV  @VVV \\
 \Spec \ k @>>> U(\varGamma)/\varGamma\ .
\end{CD}\]

Suwa [9] extended Serre's argument over commutative rings. In doing so, Serre's method, which led to the Kummer and Artin-Schreier-Witt theories by a sculpting $U(\varGamma) \rightarrow U(\varGamma)/\varGamma$, was formulated as a sculpting problem, and the embedding problem was added and formulated as follows: 

(The sculpture problem) Let $\varGamma$ be a finite group and $R$ a commutative ring. Given an affine flat $R$-group scheme $G$ and an embedding $e: \varGamma \rightarrow G$, does a commutative diagram 

\[\begin{CD}
 \varGamma   @>>> U(\varGamma) \\
 @VV{\wr}V  @VVV \\
 \varGamma @>>> G\ 
\end{CD}\]exist?

(The embedding problem) Let $\varGamma$ be a finite group and $R$ a commutative ring. Given an affine flat $R$-group scheme $G$ and an embedding $e: \varGamma \rightarrow G$, does a commutative diagram

\[\begin{CD}
 \varGamma   @>>> G \\
 @VV{\wr}V  @VVV \\
 \varGamma @>>>  U(\varGamma) \
\end{CD}\]exist?

Let $S/R$ be an unramified $\varGamma$-extension. If both the sculpture and embedding problems are affirmatively solved for $e: \varGamma \rightarrow G$, then $S/R$ has a normal basis if and only if there exist morphisms $\Spec \ S \rightarrow G$ and $\Spec \ R \rightarrow G/\varGamma$ such that the diagram 

\[\begin{CD}
 \Spec \ S  @>>> G \\
 @VVV  @VVV \\
 \Spec \ R @>>> G/\varGamma\ 
\end{CD}\]
is Cartesian.

In [9, 10], Suwa examined the problems in the case where $\varGamma$ is a cyclic group:\\
(1) Kummer theory;\\
(2) Kummer-Artin-Schreier theory; \\
(3) Artin-Schreier-Witt theory in characteristic $p>0$;\\
(4) Quadratic-twisted Kummer theory of odd degree;\\
(5) Quadratic-twisted Kummer theory of even degree;\\
(6) Quadratic-twisted Kummer-Artin-Schreier theory.
 
\vspace{3mm}

Let $S$ be a scheme and $\varGamma$ an affine commutative $S$-group scheme such that ${\mathcal O}_{\varGamma}$ is a locally free ${\mathcal O}_{S}$-module of finite rank. 

By the Cartier duality, Grothendieck devised the embedding

\[\varGamma=\mathit{Hom}_{S-\gr}(\varGamma^{\vee},\G_{m,S}) \longrightarrow \prod_{\varGamma^{\vee}/S}\G_{m,\varGamma^{\vee}}=\mathit{Hom}_{S}(\varGamma^{\vee},\G_{m,S})\]
(cf. [7, Sec 6]).
Here, $\varGamma^{\vee}$ is the Cartier dual of $\varGamma$ and  $\prod_{\varGamma^{\vee}/S}\G_{m,\varGamma^{\vee}}$ is the Weil restriction of $\G_{m,\varGamma^{\vee} }$ to $S$.
Moreover, the Weil restriction $\prod_{\varGamma^{\vee}/S}\G_{m,\varGamma^{\vee}}$ coincides with $U(\varGamma)$. Here, $U(\varGamma)$ is the unit group scheme of the group algebra of $\varGamma$. (The detailed description of $U(\varGamma)$ is recalled in the Preliminary section.)
If $\varGamma$ is a finite commutative group, the embedding invented by Serre $i: \varGamma \rightarrow U(\varGamma)$ is nothing but the above embedding $\varGamma \rightarrow \prod_{\varGamma^{\vee}/S}\G_{m,\varGamma^{\vee}}$.
 Furthermore, we have the following exact sequence of commutative $S$-group schemes
\[0 \longrightarrow \varGamma \longrightarrow U(\varGamma) \longrightarrow U(\varGamma)/\varGamma \longrightarrow 0.\]
This exact sequence is named the Grothendieck resolution.

From the above discussion, it is natural to formulate Serre's argument for a finite commutative group scheme (resp. finite non-commutative group scheme). To this end, developing a normal basis property for torsors under a group scheme was necessary. Indeed, Kreimer and Takeuchi [6] formulated the normal basis property in the framework of the Hopf-Galois theory. Moreover, Doi and Takeuchi [3] reported that the notion of Hopf-Galois extensions with the normal basis property is equivalent to the idea of cleft extensions, as defined by Sweedler [12]. Paraphrasing their arguments, we may define a cleft torsor as follows:

\vspace{3mm}

\begin{definition}[{\cite[Definition 2.9]{Tsu2}}]
Let $S$ be a scheme, $\varGamma$ an affine group $S$-scheme, and $X$ a right $\varGamma$-torsor over $S$. Then, we say that a right $\varGamma$-torsor $X$ is cleft if there exists an isomorphism of ${\mathcal O}_S$-modules $\varphi:{\mathcal O}_{\varGamma}\overset{\sim}\rightarrow {\mathcal O}_X$ such that the diagram
\[\begin{CD}
 {\mathcal O}_{\varGamma} @>{\varphi}>> {\mathcal O}_X \\
 @VV\Delta V  @VV\rho V \\
 {\mathcal O}_{\varGamma}\otimes_{{\mathcal O}_S} {\mathcal O}_{\varGamma} @>{\varphi \otimes Id}>> {\mathcal O}_X \otimes_{{\mathcal O}_S} {\mathcal O}_{\varGamma}\ .
\end{CD}\]
is commutative. Here $\Delta$ denotes the comultiplication of ${\mathcal O}_S$-Hopf algebra ${\mathcal O}_{\varGamma}$ and $\rho$ the right ${\mathcal O}_{\varGamma}$-comodule algebra structure homomorphism of ${\mathcal O}_S$-algebra ${\mathcal O}_X$. 
\end{definition}
By this definition, a characterization of cleft torsors under a finite flat group scheme is given as follows:

\vspace{3mm}
\begin{theorem}[{\cite{Su3}}]
Let $S$ be a scheme and $\varGamma$ an affine $S$-group scheme such that  ${\mathcal O}_{\varGamma}$ is a locally free ${\mathcal O}_{S}$-module of finite rank. Then, a $\varGamma$-torsor $X$ over $S$ is cleft if and only if there exists a Cartesian diagram

\[\begin{CD}
 X  @>>> U(\varGamma)\\
 @VVV  @VVV \\
 S @>>> U(\varGamma)/\varGamma\ .
\end{CD}\]
\end{theorem}
This theorem is provided by the author of this paper for the case where  $\varGamma$ is commutative, and by Suwa for the general case where  $\varGamma$ is not necessarily commutative.

\vspace{3mm}
\begin{corollary}[{\cite{Tsu2}}]Under the notation of Theorem 1.2, let $G$ be a flat affine group scheme over $S$. Suppose there exist commutative diagrams
\[\begin{CD}
 \varGamma @>i>> U(\varGamma) \\
 @VV{\wr}V  @VVV \\
 \varGamma @>>> G
\end{CD}\] 
and
\[\begin{CD}
 \varGamma @>>> G \\
 @VV{\wr}V  @VVV \\
 \varGamma @>i>> U(\varGamma)\ ,
\end{CD}\]
where $\varGamma \rightarrow G$ is a closed embedding of group schemes. Then, a $\varGamma$-torsor $X$ over $S$ is cleft if and only if $X$ is defined by a Cartesian diagram
\[\begin{CD}
 X @>>> G \\
 @VVV  @VVV \\
 S @>>> G/\varGamma\ .
\end{CD}\]
\end{corollary}

\vspace{3mm}
Let $R$ be an $\F_{p}$-algebra, $\lambda \in R$, and $W_{n,R}$ the group scheme of Witt vectors of length $n$ over $R$. Then, the isogeny $F-[\lambda^{p-1}]: W_{n, R} \rightarrow W_{n, R}$ is finite and flat. Here,  $[\lambda^{p-1}]=(\lambda^{p-1},0, \dots, 0)$ denotes the Teichm$\ddot{\text{u}}$ller representative of $\lambda^{p-1}$ in $W_{n}(R)$. 

Define $N=\Ker[F-[\lambda^{p-1}]: W_{n, R} \rightarrow W_{n, R}]$. 
In [10, Remark 4.9], Suwa demonstrated the existence of the following commutative diagrams of group schemes over $R$ with exact rows

\[\begin{CD}
 0 @>>> N @>>> U(N) @>>> U(N)/N @>>> 0 \\
 @.  @\vert  @VVV  @VVV  @. \\
 0 @>>> N @>>> W_{n,R} @>{F-[\lambda^{p-1}]}>> W_{n,R} @>>> 0
\end{CD}\]

and

\[\begin{CD}
 0 @>>> N @>>> W_{n,R} @>{F-[\lambda^{p-1}]}>> W_{n,R} @>>> 0 \\
 @.  @\vert  @VVV  @VVV  @. \\
 0 @>>> N @>>> U(N) @>>> U(N)/N @>>> 0
\end{CD}\]

Additionally, Suwa described the commutative group scheme $U(N)/N$.

\vspace{3mm}
\begin{definition} (=Definition 2.6) Let $R$ be an  $\F_{p}$-algebra. For $\lambda \in R$, we define the finite flat commutative group scheme ${\varGamma}^{(\lambda)}_{R}$ over $R$ by  ${\varGamma}^{(\lambda)}_{R}= \Spec R[T]/(T^{p^{n}})$ with
\begin{align*}
&\text{(a) the multiplication}: T \mapsto T \otimes 1+1 \otimes T+\lambda T \otimes T; \ \ \ \ \ \ \ \ \ \ \ \ \ \ \ \ \ \ \ \ \ \ \ \ \ \ \\
&\text{(b) the unit}: T \mapsto 0;  \\
&\text{(c) the inverse}: T \mapsto -\frac{1}{1+\lambda T}.
\end{align*}
\end{definition}
In this paper, the $R$-Hopf algebra $R[T]/(T^{p^{n}})$ representing ${\varGamma}^{(\lambda)}_{R}$ is denoted by $A_{R}^{(\lambda)}$. 

Moreover, the Cartier dual of ${\varGamma}^{(\lambda)}_{R}$ is $N$. This was demonstrated by the author in [14] for the case of $n=1$ and by Amano [1] for the general case where $n$ is any natural number.

The main results of this paper are as follows:

\vspace{5mm}

\begin{theorem} (=Theorem 3.8) We have commutative diagrams of group schemes over $R$ with exact rows:

\[\begin{CD}
 0 @>>> {\varGamma}^{(\lambda)}_{R} @>e>>  {\mathcal G}^{(\lambda)}_{R}   @>{F^{n}}>> \mathcal{G}_{R}^{(\lambda^{p^{n}})} @>>> 0\\
 @.     @\vert     @VV\sigma_{1}V    @VV\tau_{1}V  \\
 0 @>>> {\varGamma}^{(\lambda)}_{R} @>i>>  U({\varGamma}^{(\lambda)}_{R})  @>q>>  U({\varGamma}^{(\lambda)}_{R}) / {\varGamma}^{(\lambda)}_{R} @>>> 0
\end{CD}\]
and
\[\begin{CD}
 0 @>>> {\varGamma}^{(\lambda)}_{R} @>i>>  U({\varGamma}^{(\lambda)}_{R})   @>q>> U({\varGamma}^{(\lambda)}_{R}) / {\varGamma}^{(\lambda)}_{R} @>>> 0\\
 @.     @\vert     @VV\sigma_{2}V    @VV\tau_{2}V  \\
 0 @>>> {\varGamma}^{(\lambda)}_{R} @>e>>  {\mathcal G}^{(\lambda)}_{R} @>{F^{n}}>> {\mathcal G}^{(\lambda^{p^{n}})}_{R} @>>> 0.
\end{CD}\]

\end{theorem}

Additionally, the commutative group scheme $U({\varGamma}^{(\lambda)}_{R}) / {\varGamma}^{(\lambda)}_{R}$ is described in Theorem 3.7. Kassel and Masuoka [5] addressed the Noether problem for Hopf algebras, with Theorem 3.7 providing a noteworthy positive example. 
We conclude this article by presenting an example of non-cleft  ${\varGamma}^{(\lambda)}_{R}$-torsors.

\vspace{5mm}

\noindent{\bf Notation 1.6.} In this article, $p$ denotes a prime number and $n$ denotes a positive integer. For a scheme $S$ and a group scheme $\varGamma$ over $S$, $H^{1}(S, \varGamma)$ denotes the set of isomorphism classes of right $\varGamma$-torsors over $S$. (For details, please refer to Demazure-Gabriel [2, Ch III. 4].)

\section{Preliminary}\label{sec:prelimi}

\vspace{5mm}

\subsection{$U(\varGamma)$ for a finite flat group scheme $\varGamma$.}
We recall the group algebra scheme $A(\varGamma)$ and its unit group scheme $U(\varGamma)$ for a finite flat group scheme $\varGamma$. For details of these group schemes, we refer to [11, Section 2].

Let $S$ be a scheme and $\varGamma$ an affine group scheme over $S$. 
Let $A(\varGamma)$ denote the ring functor defined by $T \mapsto \mathit{Hom}_{{\mathcal{O}}_{S}}({\mathcal{O}}_{\varGamma}, {\mathcal{O}}_{T})$ for an affine $S$-scheme $T$, where the multiplication of $\mathit{Hom}_{{\mathcal{O}}_{S}}({\mathcal{O}}_{\varGamma}, {\mathcal{O}}_{T})$ is defined by the convolution product. Then, $A(\varGamma)$ is an $S$-ring scheme. Moreover, we define a functor $U(\varGamma)$ by $U(\varGamma)(T)=A(\varGamma)(T)^{\times}$ for an affine $S$-scheme $T$. Then, $U(\varGamma)$ is a sheaf of groups for the fppf-topology over $S$. If ${\mathcal{O}}_{\varGamma}$ is a locally free ${\mathcal{O}}_{S}$-module of finite rank, $U(\varGamma)$ is represented by an affine smooth group scheme over $S$.
Let $R$ be a commutative ring. We assume that $S=\Spec\,R$ and $\varGamma=\Spec\,H$, where $H$ is a free $R$-module of finite rank. We take a basis $\{e_{1}, \dots, e_{n}\}$ of $H$ over $R$. Let $S_{R}(H)$ denote the symmetric $R$-algebra associated with the $R$-module $H$. For each $i$, let $T_{e_{i}}$ denote the image of $e_{i}$ by the canonical injection $H \rightarrow S_{R}(H)$. Moreover, we define a linear combination $R_{ij}(e_{1}, \dots, e_{n})=\displaystyle \sum_{k=1}^{n}c_{ijk}e_{k}$ for each $1 \leq i, j \leq n$ by
\[\varDelta_{H}(e_{j})=\displaystyle \sum_{i=1}^{n}e_{i} \otimes R_{ij}(e_{1}, \dots , e_{n}).\]
Then, we obtain that $A(\varGamma)=\Spec\,S_{R}(H)=\Spec\,R[T_{e_{1}}, \dots T_{e_{n}}]$ with
the comultiplication map
\[\varDelta(T_{e_{j}})=\displaystyle \sum_{i=1}^{n}T_{e_{i}} \otimes R_{ij}(T_{e_{1}}, \dots, T_{e_{n}}),\]
where $R_{ij}(T_{e_{1}}, \dots, T_{e_{n}})=\displaystyle \sum_{k=1}^{n}c_{ijk}T_{e_{k}}$ and the counit map $\varepsilon(T_{e_{j}})=\varepsilon_{H}(e_{j})$. Moreover, let $A$ be an $R$-algebra. Then, the multiplication of $A(\varGamma)(A)$ is defined by
\[(a_{1}, a_{2}, \dots, a_{n})(b_{1}, b_{2}, \dots, b_{n})\]
\[=(\displaystyle \sum_{j=1}^{n}R_{1j}(a_{1},a_{2}, \dots, a_{n})b_{j},  \sum_{j=1}^{n}R_{2j}(a_{1},a_{2}, \dots, a_{n})b_{j}, \dots,\sum_{j=1}^{n}R_{nj}(a_{1},a_{2}, \dots, a_{n})b_{j} )\]
Hence, $(a_{1}, a_{2}, \dots, a_{n}) \in A(\varGamma)(A)$ is invertible if and only if $\det(R_{ij}(a_{1}, a_{2}, \dots, a_{n}))$ is invertible in $A$.
Therefore, we have 
\[U(\varGamma)=\Spec\,R[T_{e_{1}}, T_{e_{2}}, \dots T_{e_{n}}, \frac{1}{D}],\]
where $D=\det(R_{ij}(T_{e_{1}}, T_{e_{2}}, \dots, T_{e_{n}}))$.

Moreover, the $R$-homomorphism $i^{\#}:  R[T_{e_{1}}, T_{e_{2}}, \dots T_{e_{n}}, {1}/{D}] \rightarrow H$ defined by $T_{e_{i}} \mapsto e_{i}$ induces a closed immersion of group schemes
$i: \varGamma \rightarrow U(\varGamma)$.  

The above-mentioned result is the sketch of the proof of [11, Theorem 2.6]. If $R$ is a field, the above Hopf algebra $R[T_{1}, T_{2}, \dots T_{n}, {1}/{D}]$ coincides the commutative free Hopf algebra generated by $H$ constructed by Takeuchi [13]. Therefore, [11] gives an algebraic geometric interpretation of Takeuchi's result.

\vspace{5mm}

\vspace{3mm}

\subsection{The group schemes considered in this article.} In this section, $R$ denotes an $\F_{p}$-algebra.

\vspace{3mm}

\begin{definition} The additive group scheme $\G_{a, R}$ over $R$ is defined by
\[\G_{a, R}=\Spec\ R[T]\] 
with
\begin{align*}
&\text{(a) the multiplication}: T \mapsto T\otimes 1 + 1 \otimes T; \ \ \ \ \ \ \ \ \ \ \ \ \ \ \ \ \ \ \ \ \ \ \ \ \ \ \ \ \ \ \ \ \ \ \ \ \ \ \ \ \ \ \ \ \ \ \ \ \ \ \ \ \ \ \ \ \\
&\text{(b) the unit}: T \mapsto 0;\\
&\text{(c) the inverse}: T \mapsto -T.
\end{align*}

\end{definition}
\vspace{3mm}

\begin{definition} The multiplicative group scheme $\G_{m, R}$ over $R$ is defined by
\[\G_{m, R}=\Spec\ R[U, \frac{1}{U}]\] 
with
\begin{align*}
&\text{(a) the multiplication}: U \mapsto U\otimes U; \ \ \ \ \ \ \ \ \ \ \ \ \ \ \ \ \ \ \ \ \ \ \ \ \ \ \ \ \ \ \ \ \ \ \ \ \ \ \ \ \ \ \ \ \ \ \ \ \ \ \ \ \ \ \ \ \ \ \ \ \ \ \ \\
&\text{(b) the unit}: U \mapsto 1;\\
&\text{(c) the inverse}: U \mapsto \frac{1}{U}.
\end{align*}

\end{definition}
\vspace{3mm}

\begin{definition} The finite flat commutative group scheme ${\boldsymbol{\boldsymbol{\mu}}_{p^{n}, R}}$ over $R$ is defined by  \[\boldsymbol{\mu}_{{p^{n}, R}} = \Spec R[U]/(U^{p^{n}}-1)\] with
\begin{align*}
&\text{(a) the multiplication}: U \mapsto U\otimes U; \ \ \ \ \ \ \ \ \ \ \ \ \ \ \ \ \ \ \ \ \ \ \ \ \ \ \ \ \ \ \ \ \ \ \ \ \ \ \ \ \ \ \ \ \ \ \ \ \ \ \ \ \ \ \ \ \ \ \ \ \ \ \ \\
&\text{(b) the unit}: U \mapsto 1;\\
&\text{(c) the inverse}: U \mapsto \frac{1}{U}.
\end{align*}
\end{definition}
\vspace{3mm}

In this paper, the $R$-Hopf algebra $R[U]/(U^{p^{n}}-1)$ representing ${{\boldsymbol{\mu}_{p^{n}}}_{R}}$ is denoted by $A_{\boldsymbol{\mu}, R}$. 

Given the $p^{n}$-power map $F^{n}: \G_{m, R} \rightarrow  \G_{m, R}$,  we obtain that ${\boldsymbol{\mu}_{p^{n}, R}}=\Ker [F^{n}: \G_{m, R} \rightarrow  \G_{m, R}]$. 

Therefore, we have an exact sequence of $R$-group schemes

\[0 \rightarrow {\boldsymbol{\mu}_{p^{n}, R}} \rightarrow \G_{m, R}  \overset{F^{n}}\rightarrow \G_{m, R} \rightarrow 0.\]
This sequence is the famous Kummer sequence.  

\vspace{3mm}

\begin{definition} For $\lambda \in R$, the commutative group scheme $\mathcal{G}_{R}^{(\lambda)}$ over $R$ is defined by

\[\mathcal{G}_{R}^{(\lambda)} = \Spec\ R[T, \frac{1}{1+\lambda T}]\]
with
\begin{align*}
&\text{(a) the multiplication}: T \mapsto T\otimes 1 +1 \otimes T+ \lambda T \otimes T; \ \ \ \ \ \ \ \ \ \ \ \ \ \ \ \ \ \ \ \ \ \ \ \ \ \ \ \ \ \ \ \ \ \ \ \ \ \ \ \ \ \ \ \ \ \\
&\text{(b) the unit}: T \mapsto 0;\\
&\text{(c) the inverse}: T \mapsto -\frac{T}{1+\lambda T}.
\end{align*}

\end{definition}
\vspace{3mm}

A homomorphism $\alpha^{(\lambda)}: \mathcal{G}_{R}^{(\lambda)} \rightarrow \G_{m, R}$ of group schemes over $R$ is defined by
\[U \mapsto \lambda T+1 : R[U, \frac{1}{U}] \rightarrow R[T, \frac{1}{1+\lambda T}]\]
If $\lambda$ is invertible in $R$, then $\alpha^{(\lambda)}$ is an isomorphism. Conversely, if $\lambda=0, \mathcal{G}_{R}^{(\lambda)}$ is nothing but the additive group scheme   $\G_{a, R}$.     

\vspace{3mm}
\begin{definition} For $\lambda \in R$, we define the finite flat commutative group scheme $\varGamma^{(\lambda)}_{R}$ over $R$ by  $\varGamma^{(\lambda)}_{R}= \Spec R[T]/(T^{p^{n}})$ with
\begin{align*}
&\text{(a) the multiplication}: T \mapsto T\otimes 1 +1 \otimes T+ \lambda T \otimes T; \ \ \ \ \ \ \ \ \ \ \ \ \ \ \ \ \ \ \ \ \ \ \ \ \ \ \ \ \ \ \ \ \ \ \ \ \ \ \ \ \ \ \ \ \ \\
&\text{(b) the unit}: T \mapsto 0;\\
&\text{(c) the inverse}: T \mapsto -\frac{T}{1+\lambda T}.
\end{align*}
\end{definition}
In this paper, the $R$-Hopf algebra $R[T]/(T^{p^{n}})$ representing $\varGamma^{(\lambda)}_{R}$ is denoted by $A_{R}^{(\lambda)}$.

\vspace{3mm}
\subsection{}
Given $p^{n}$-power map $F^{n}: \mathcal{G}_{R}^{(\lambda)} \rightarrow  \mathcal{G}_{R}^{(\lambda^{p^{n}})}$, we have that $\varGamma^{(\lambda)}_{R}=\Ker [F^{n}: \mathcal{G}_{R}^{(\lambda)} \rightarrow  \mathcal{G}_{R}^{(\lambda^{p^{n}})}]$. 

Therefore, we obtain an exact sequence of $R$-group schemes 

\[0 \rightarrow \varGamma^{(\lambda)}_{R} \rightarrow  \mathcal{G}_{R}^{(\lambda)} \overset{F^{n}}\rightarrow  \mathcal{G}_{R}^{(\lambda^{p^{n}})} \rightarrow 0.\]

\vspace{3mm}

First, we describe the unit group scheme $U({\boldsymbol{\mu}_{p^{n}, R}})$ of the group ring of ${\boldsymbol{\mu}_{p^{n}, R}}$.

\vspace{3mm}

\begin{prop} The group algebra scheme for the finite flat group scheme ${\boldsymbol{\mu}_{p^{n}, R}}$ is given by
\[A({\boldsymbol{\mu}_{p^{n}, R}})=\Spec\,R[Y_{1}, Y_{U}, Y_{U^{2}}, \dots Y_{U^{p^{n}-1}}]\]
with

(a) the multiplication:
\[Y_{1} \mapsto Y_{1} \otimes Y_{1},\]
\[Y_{U^{s}} \mapsto Y_{U^{s}} \otimes Y_{U^{s}} \  \text{for} \  1 \leq s \leq p^{n}-1;\]

(b) the unit:
\[Y_{1} \mapsto 1,\]
\[Y_{U^{s}} \mapsto 1  \  \text{for} \  1 \leq s \leq p^{n}-1. \]

Moreover, the unit group scheme of the group algebra for the finite flat group scheme ${\boldsymbol{\mu}_{p^{n}, R}}$ is given by
\[U({\boldsymbol{\mu}_{p^{n}, R}})=\Spec\,R[Y_{1}, Y_{U}, Y_{U^{2}}, \dots Y_{U^{p^{n}-1}}, \frac{1}{{D}'}] ,\]
where

 \[{D}'=Y_{1} Y_{U} Y_{U^{2}} \dots Y_{U^{p^{n}-1}}. \]

\end{prop}
\begin{proof} Note that 
\[\beta=\{1, U^{1}, U^{2}, \dots, U^{p^{n}-1}\}\] is a basis of $R$-module  $A_{\boldsymbol{\mu}, R}$. 

We have 
\begin{align*}
\Delta(U^{s})&=U^{s} \otimes U^{s},\\
\end{align*}
where $\Delta$ is the multiplication map of  $A_{\boldsymbol{\mu}, R}$.

We have 
\begin{align*}
\varepsilon(Y^{s})&=1,\\
\end{align*}
where $\varepsilon$ is the unit of $A_{\boldsymbol{\mu}, R}$.

Considering $e_{i+1}=U^{i}$ for $0 \leq i \leq p^{n}-1$, from Subsection 2.1,
we know the multiplication and unit of $A(\boldsymbol{\mu}_{p^{n}, R})$.
Furthermore, the right regular representation $(R_{ij})_{1 \leq i,j \leq p^{n}}$ of $A_{\boldsymbol{\mu}, R}$ with respect to the basis $\beta$ is given by 

\[R_{ij}(U^{0}, U^{1}, U^{2}, \dots, U^{p^{n}-1})=\begin{cases}
0 & (i>j)\\
U^{i-1} & (i=j)\\ 
0 & (i<j).
\end{cases}\] 
Therefore, we have 
\[{D}'=\det\Big(R_{ij}(Y_{U^{0}}, Y_{U^{1}}, \dots ,Y_{U^{p^{n}-1}})\Big)=Y_{U^{0}}Y_{U^{1}}\dots Y_{U^{p^{n}-1}}\]
because the matrix $(R_{ij})_{1 \leq i,j \leq p^{n}}$ is a diagonal matrix.

\end{proof}

\vspace{3mm}
In this paper, the $R$-Hopf algebra $R[Y_{1}, Y_{U}, Y_{U^{2}}, \dots Y_{U^{p^{n}-1}}, {1}/{D}']$ representing $U(\boldsymbol{\mu}_{p^{n}, R})$ is denoted by $S(A_{\boldsymbol{\mu}, R})_{\Theta}$. 
\vspace{3mm}

Next, we describe the unit group scheme $U(\varGamma^{(\lambda)}_{R})$ of the group ring of $\varGamma^{(\lambda)}_{R}$.

\vspace{3mm}

\begin{prop} For $\lambda \in R$, the group algebra scheme for $\varGamma^{(\lambda)}_{R}$ is defined by
\[A(\varGamma^{(\lambda)}_{R})=\Spec\, R[X_{1}, X_{T}, X_{T^{2}}, \dots, X_{T^{p^{n}-1}}]\] 
with

(a) the multiplication:

\[X_{1} \mapsto X_{1} \otimes X_{1},\]
\[X_{T^{s}} \mapsto \sum_{l=0}^{s}\binom{s}{l}\Big\{\sum_{k=0}^{l}{\binom{l}{k}}\lambda^{k}X_{T^{k+s-l}}\Big\}\otimes X_{T^{l}} \  \text{for} \  1 \leq s \leq p^{n}-1;\]

(b) the unit:

\[X_{1} \mapsto 1,\]
\[X_{T^{s}} \mapsto 0  \  \text{for} \  1 \leq s \leq p^{n}-1. \]

Moreover, the unit group scheme of the group algebra for the finite flat group scheme $\varGamma^{(\lambda)}$ is defined by
\[U(\varGamma^{(\lambda)}_{R})=\Spec\ R[X_{1}, X_{T}, X_{T^{2}}, \dots, X_{T^{p^{n}-1}}, \frac{1}{D_{(\lambda)}}],\]
where \[D_{(\lambda)}=X_{1}\prod_{l=1}^{p^{n}-1} \Big\{\sum_{k=0}^{l}{\binom{l}{k}}\lambda^{k}X_{T_{k}}\Big\}.  \]
\end{prop}
\vspace{3mm}

\begin{proof} Note that 
\[\beta^{(\lambda)}=\{1, T^{1}, T^{2}, \dots, T^{p^{n}-1}\}\] is a basis of $R$-module  $A^{(\lambda)}_{R}$. 

We have 
\begin{align*}
\Delta^{(\lambda)}(T^{0})& = T^{0} \otimes T^{0},\\
\Delta^{(\lambda)}(T^{s}) &= \sum_{l=0}^{s}{\binom{s}{l}}\Big\{\sum_{k=0}^{l}{\binom{l}{k}}\lambda^{k}{T^{k+s-l}}\Big\}\otimes T^{l} \  \text{for} \  1 \leq s \leq p^{n}-1,\\
\end{align*}
where $\Delta^{(\lambda)}$ is the multiplication of $A_{R}^{(\lambda)}$.

We have 
\begin{align*}
\varepsilon^{(\lambda)}(T^{0})&=1,\\
\varepsilon^{(\lambda)}(T^{s})&=0  \  \text{for} \  1 \leq s \leq p^{n}-1,\\
\end{align*}
where $\varepsilon^{(\lambda)}$ is the counit map of $A^{(\lambda)}_{R}$.

Considering $e_{i+1}=T^{i}$, for $0 \leq i \leq p^{n}-1$, from Subsection 2.1, we obtain the multiplication and unit of $A(\varGamma^{(\lambda)}_{R})$.
Furthermore, the right regular representation $(R_{ij})_{1 \leq i,j \leq p^{n}}$ of $A^{(\lambda)}_{R}$ with respect to the basis $\beta^{(\lambda)}$ is given by
\[R_{ij}(T^{0}, \dots, T^{p^{n}-1})=\begin{cases}
0 & (i>j). \\[2mm]
\displaystyle\sum_{k=0}^{i-1}{\binom{i-1}{k}}\lambda^{k}T^{k} & (i=j)\\
\text{a polynomial of $T$} & (i<j)
\end{cases}\]

Therefore, we have that  
\[D_{(\lambda)}=R_{ij}(X_{T^{0}}, \dots, X_{T^{p^{n}-1}})=\displaystyle\prod_{r=0}^{p^{n}-1}\Big\{\displaystyle\sum_{k=0}^{r}{\binom{r}{k}}\lambda^{k}X_{T^{k}}\Big\}\]
since $R_{ij}(T^{0}, \dots, T^{p^{n}-1})$ is an upper triangular matrix.
\end{proof}
\vspace{5mm}

In this paper, the $R$-Hopf algebra $R[X_{1}, X_{T}, X_{T^{2}}, \dots, X_{T^{p^{n}-1}}, {1}/{D_{(\lambda)}}]$ representing $U(\varGamma^{(\lambda)}_{R})$ is denoted by $S(A_{R}^{(\lambda)})_{\Theta}$. 

\vspace{3mm}

The homomorphism $\alpha^{(\lambda)}: \mathcal{G}_{R}^{(\lambda)} \rightarrow \G_{m, R}$ of group schemes over $R$ induces $\widetilde{\alpha^{(\lambda)}}: \varGamma^{(\lambda)}_R \rightarrow {\boldsymbol{\mu}_{p^{n}, R}}$. This is defined by 
\[R[U]/(U^{p^{n}}-1) \rightarrow R[T]/(T^{p^{n}}), \]
\[U \mapsto 1+\lambda T.\]

If $\lambda$ is invertible in $R$, then $\widetilde{\alpha^{(\lambda)}}$ is an isomorphism.

Moreover, $\widetilde{\alpha^{(\lambda)}}$ induces $U(\widetilde{\alpha^{(\lambda)}}): U(\varGamma^{(\lambda)}_R) \rightarrow U( {\boldsymbol{\mu}_{p^{n}, R}})$. This is defined by
\[Y_{1} \mapsto X_{1},\]
\[Y_{U^{s}} \mapsto \sum_{k=0}^{s}{\binom{s}{k}}\lambda^{k}X_{T^{s}} \ \text{for} \ 1 \leq s \leq p^{n}-1.\]

\vspace{3mm}
\section{Main result}\label{sec:main}

\vspace{3mm}

To prove the main results, we recall the results of Doi and Takeuchi. Let $R$ be a commutative ring, $H$ an $R$-Hopf algebra, and $A$ a right $H$-comodule $R$-algebra (with $H$-comodule algebra structure $\rho: A \rightarrow A \otimes H$).
We denote by $A^{\mathrm{co}H}=\{a \in A | \rho(a)=a \otimes 1\}$ the coinvariant subalgebra in $A$ by the right $H$-coaction of $A$.
If there exists a homomorphism of $R$-modules $\phi: H \rightarrow A$, which is also a homomorphism of right $H$-comodules and invertible for the convolution product, the extension $A/A^{\mathrm{co}H}$ of $R$-algebras is called a cleft extension ([12]). Moreover, the homomorphism of $R$-modules $\phi$ is called a cleaving map of $A$. 
Then, a left $A^{\mathrm{co}H}$-module homomorphism $\Phi : A^{\mathrm{co}H} \otimes H \rightarrow A$ defined by $b \otimes h \mapsto b\phi(h) $ is bijective. Moreover, we define a homomorphism of $R$-modules $P : A \rightarrow A$ by $a \mapsto \sum a_{(0)}\phi^{-1}(a_{(1)})$ for $a \in A.$ Then, we have $P(A) \subset A^{\mathrm{co}H}$. Furthermore, \[\Phi^{-1}(a)=\sum P(a_{(0)})\otimes a_{(1)}  \ \ \ \ \ \ \ \ \ \   (*)\] (cf. [3, Theorem 9]).

\vspace{3mm}
The following results are important in deriving the main results.
\vspace{3mm}

Let $S$ be a scheme and $\varGamma$ an affine commutative group $S$-scheme such that $\mathcal{O}_{\varGamma}$ is a locally free of $\mathcal{O}_{S}$-module of finite rank. Then, by the natural closed immersion $i: \varGamma \rightarrow U(\varGamma)$, $U(\varGamma)$ is a cleft $\varGamma$-torsor over $U(\varGamma)/\varGamma$ (see [11, Proposition 3.1]).

If $S= \Spec \ R$ and $\varGamma= \Spec \ H$, where $H$ is a free $R$-module of finite rank. Then, by taking a basis$\{e_{1}, \dots e_{n}\}$ of $H$ over $R$, the canonical closed immersion $i: \varGamma \rightarrow U(\varGamma)$ is given by

\[i^{\#} : R[T_{e_{1}}, \dots, T_{e_{n}}, \frac{1}{D}]\rightarrow H,\]
\[T_{e_{i}} \mapsto e_{i} \ \ \text{for} \ \ 1 \leq i \leq n.\]
By the $R$-homomorphism $i^{\#}$ that determines the canonical closed immersion $i$, ${S(H)_{\Theta}}$ becomes a right $H$-comodule algebra. Furthermore, a homomorphism of $R$-modules \[\phi: H \rightarrow  {S(H)_{\Theta}}\]
\[e_{i} \mapsto T_{e_{i}}  \ \ \text{for} \ \ 1 \leq i \leq n.\]
is a cleaving map of ${S(H)_{\Theta}}$.
Moreover, we see that $U(\varGamma)/\varGamma$ is isomorphic to $\Spec \ {S(H)_{\Theta}}^{\mathrm{co}H}$ as $R$-commutative group schemes since $H$ is finite over $R$. The canonical morphism $q: U(\varGamma) \rightarrow U(\varGamma)/\varGamma$ is given by the injection ${S(H)_{\Theta}}^{\mathrm{co}H} \subset {S(H)_{\Theta}}$ (cf. [4, Part I, 5.6, 6.1]) .

\vspace{3mm}

In this section, $R$ denotes an $\F_{p}$-algebra.

\vspace{3mm}

For example, we consider for $\boldsymbol{\mu}_{p^{n}, R}$.
The canonical closed immersion $i: \boldsymbol{\mu}_{p^{n}, R} \rightarrow U({\boldsymbol{\mu}_{p^{n}, R}})$ is given by
\[ i^{\#}: R[Y_{1}, Y_{U}, Y_{U^{2}}, \dots Y_{U^{p^{n}-1}}, \frac{1}{{D}'}  ] \rightarrow R[U]/(U^{p^{n}}-1),\]
\[Y_{1} \mapsto 1,  \  Y_{U^{s}} \mapsto U^{s}  \  \text{for} \  1 \leq s \leq p^{n}-1. \]
By the $R$-Hopf algebras map $i^{\#}$ that determines this canonical closed immersion $i$, $S(A_{\boldsymbol{\mu}, R})_{\Theta}$ becomes a right $A_{\boldsymbol{\mu}, R}$-comodule algebra.

Moreover, a homomorphism of $R$-modules \[\phi: R[U]/(U^{p^{n}}-1) \rightarrow R[Y_{1}, Y_{U}, Y_{U^{2}}, \dots Y_{U^{p^{n}-1}}, \frac{1}{{D}'}, ]\]
\[U^{s} \mapsto Y_{U^{s}} \  \  0 \leq s \leq p^{n}-1\]
is a right $A_{\boldsymbol{\mu}, R}$-comodule map and is invertible for the convolution product.

Therefore, a homomorphism of left ${S(A_{\boldsymbol{\mu}, R})_{\Theta}}^{\mathrm{co}A_{\boldsymbol{\mu}, R}}$-modules, which is also a right $A_{\boldsymbol{\mu}, R}$-comodule map,
\[\Phi:{S(A_{\boldsymbol{\mu}, R})_{\Theta}}^{\mathrm{co}A_{\boldsymbol{\mu}, R}}  \otimes A_{\boldsymbol{\mu}, R} \rightarrow S(A_{\boldsymbol{\mu}, R})_{\Theta}, \]
\[b \otimes a \mapsto b\phi(a)\]
is bijective.

\vspace{3mm}

Moreover, we consider $\varGamma^{(\lambda)}_R$. The canonical closed immersion $j: \varGamma^{(\lambda)}_R \rightarrow U(\varGamma^{(\lambda)}_R)$ is given by

\[ j^{\#}: R[X_{1}, X_{T}, X_{T^{2}}, \dots X_{T^{p-1}}, \frac{1}{{D_{(\lambda)}}}  ] \rightarrow R[T]/(T^{p^{n}}),\]
\[X_{1} \mapsto 1,  \  X_{T^{s}} \mapsto T^{s}  \  \text{for} \  1 \leq s \leq p^{n}-1. \]

By the homomorphism of $R$-Hopf algebras $j^{\#}$ that determines this canonical closed immersion $j$, $S(A_{R}^{(\lambda)})_{\Theta}$ becomes a right $A_{R}^{(\lambda)}$-comodule algebra. 
Moreover, a homomorphism of $R$-modules \[\psi: R[T]/(T^{p}) \rightarrow R[X_{1}, X_{T}, Y_{X^{2}}, \dots X_{T^{p^{n}-1}}, \frac{1}{{D_{(\lambda)}}} ],\]
\[T^{s} \mapsto X_{T^{s}} \  \  0 \leq s \leq p^{n}-1\]
is a right $A_{R}^{(\lambda)}$-comodule map and is invertible for the convolution product.

Therefore, a homomorphism of left $S(A_{R}^{(\lambda)})_{\Theta}^{\mathrm{co}A_{R}^{(\lambda)}}$-modules, which is also a right $A_{R}^{(\lambda)}$-comodule map,
\[\Psi: S(A_{R}^{(\lambda)})_{\Theta}^{\mathrm{co}A_{R}^{(\lambda)}} \otimes A_{R}^{(\lambda)} \rightarrow S(A_{R}^{(\lambda)})_{\Theta}, \]
\[b \otimes a \mapsto b\psi(a)\]
is bijective.

\vspace{3mm}

\begin{prop} \[{S(A_{\boldsymbol{\mu}, R})_{\Theta}}^{\mathrm{co}A_{\boldsymbol{\mu}, R}}=R[Y_{1}^{\pm 1}, (\frac{{Y_{U}}^{2}}{Y_{U^{2}}})^{\pm 1}, \dots  (\frac{{Y_{U}}^{p^{n}-1}}{Y_{U^{p^{n}-1}}})^{\pm 1}, ({Y_{U}}^{p^{n}})^{\pm 1} ].\]
\end{prop}
\vspace{3mm}

\begin{proof} First, we have that 
\[R[Y_{1}^{\pm 1}, (\frac{{Y_{U}}^{2}}{Y_{U^{2}}})^{\pm 1}, \dots  (\frac{{Y_{U}}^{p^{n}-1}}{Y_{U^{p^{n}-1}}})^{\pm 1}, ({Y_{U}}^{p^{n}})^{\pm 1} ] \subset {S(A_{\boldsymbol{\mu}, R})_{\Theta}}^{\mathrm{co}A_{\boldsymbol{\mu}, R}}. \]
Since any integer $N$ is represented as $N=p^{n}m+r$ for some integer $m$ and some integer $r$ such that $0 \leq r < p^{n}$,

\[\Phi^{-1}({Y_{U}}^{N})=\Phi^{-1}({Y_{U}}^{p^{n}m+r})=({Y_{U}}^{p^{n}})^{m}\frac{{Y_{U}}^r}{Y_{U^{r}}} \otimes U^{r}.\]

Moreover, for $m_{0}, m_{1}, \dots, m_{p^{n}-1} \in \Z$, 

\[{Y_{1}}^{m_{0}}{Y_{U}}^{m_{1}} \dots {Y_{U^{p^{n}-1}}}^{m_{p^{n}-1}}={Y_{1}}^{m_{0}}(\frac{Y_{U^{2}}}{{Y_{U}}^{2}})^{m_{2}} \dots (\frac{Y_{U^{p^{n}-1}}}{{Y_{U}}^{p^{n}-1}})^{m_{p^{n}-1}}{Y_{U}}^{m_{1}+\sum_{k=2}^{p^{n}-1}km_{k}}.\]

Hence, we obtain
\[{\Phi}^{-1}({Y_{1}}^{m_{0}}{Y_{U}}^{m_{1}} \dots {Y_{U^{p^{n}-1}}}^{m_{p^{n}-1}})\]\[=\Big\{{Y_{1}}^{m_{0}}(\frac{Y_{U^{2}}}{{Y_{U}}^{2}})^{m_{2}} \dots (\frac{Y_{U^{p^{n}-1}}}{{Y_{U}}^{p^{n}-1}})^{m_{p^{n}-1}} \otimes 1\Big\} {\Phi}^{-1}({Y_{U}}^{m_{1}+\sum_{k=2}^{p^{n}-1}km_{k}}) .\]

Therefore, a homomorphism of $R$-modules \[{\Phi}': R[Y_{1}^{\pm 1}, (\frac{{Y_{U}}^{2}}{Y_{U^{2}}})^{\pm 1}, \dots  (\frac{{Y_{U}}^{p^{n}-1}}{Y_{U^{p^{n}-1}}})^{\pm 1}, ({Y_{U}}^{p^{n}})^{\pm 1} ] \otimes A_{\boldsymbol{\mu}, R} \rightarrow S(A_{\boldsymbol{\mu}, R})_{\Theta},  \]
\[b \otimes a \mapsto b \phi(a) \]
is bijective.

Let $s$ be \[R[Y_{1}^{\pm 1}, (\frac{{Y_{U}}^{2}}{Y_{U^{2}}})^{\pm 1}, \dots (\frac{{Y_{U}}^{p^{n}-1}}{Y_{U^{p^{n}-1}}})^{\pm 1}, ({Y_{U}}^{p^{n }})^{\pm 1} ]\] to ${S(A_{\boldsymbol{\mu}, R})_{\Theta}}^{\mathrm{co}A_{\boldsymbol{\mu}, R}}$ as the natural embedding. Then, we obtain the following commutative diagram 

\[
  \xymatrix{
     R[Y_{1}^{\pm 1}, (\frac{{Y_{U}}^{2}}{Y_{U^{2}}})^{\pm 1}, \dots  (\frac{{Y_{U}}^{p^{n}-1}}{Y_{U^{p-1}}})^{\pm 1}, ({Y_{U}}^{p^{n}})^{\pm 1} ] \otimes A_{\boldsymbol{\mu}, R}\ar[r]^{{ \ \ \ \ \ \ \ \ \ \ \ \ \ \ \ \ \ \ \ \ \ \ \ \ \Phi}'} \ar[d]_{s \otimes \text{Id}} &  S(A_{\boldsymbol{\mu}, R})_{\Theta}  \\
     {S(A_{\boldsymbol{\mu}, R})_{\Theta}}^{\mathrm{co}A_{\boldsymbol{\mu}, R}}  \otimes A_{\boldsymbol{\mu}, R}  \ar[ur]^{\Phi}. & 
  }
\]

$s$ is bijecive since $A_{\boldsymbol{\mu}, R}$ is faithfully flat over $R$.

\end{proof}

\vspace{3mm}

\begin{prop}
The $R$-subalgebra \[R[Y_{1}^{\pm 1}, (\frac{{Y_{U}}^{2}}{Y_{U^{2}}})^{\pm 1}, \dots  (\frac{{Y_{U}}^{p^{n}-1}}{Y_{U^{p^{n}-1}}})^{\pm 1}, ({Y_{U}}^{p^{n}})^{\pm 1} ]\] of  $S(A_{\boldsymbol{\mu}, R})_{\Theta}$
is isomorphic to the $R$-algebra $R[{y_{1}}^{\pm 1}, {y_{2}}^{\pm 1}, \dots {y_{p^{n}}}^{\pm 1}]$.
\end{prop}

\begin{proof}
We define a homomorphism of $R$-algebras \[\xi: R[{y_{1}}^{\pm 1}, {y_{2}}^{\pm 1}, \dots {y_{p^{n}}}^{\pm 1}] \rightarrow R[Y_{1}^{\pm 1}, (\frac{{Y_{U}}^{2}}{Y_{U^{2}}})^{\pm 1}, \dots  (\frac{{Y_{U}}^{p^{n}-1}}{Y_{U^{p^{n}-1}}})^{\pm 1}, ({Y_{U}}^{p^{n}})^{\pm 1} ]\]
by
\[y_{1} \mapsto Y_{1}, \ y_{2} \mapsto \frac{{Y_{U}}^{2}}{Y_{U^{2}}}, \dots, y_{p^{n}-1} \mapsto \frac{{Y_{U}}^{p^{n}-1}}{Y_{U^{p^{n}-1}}}, \ y_{p^{n}} \mapsto {Y_{U}}^{p^{n}}. \] 

This is well-defined and, additionally, we know that $\xi$ is bijective.

\end{proof}
\vspace{3mm}
Moreover, the $R$-algebra $R[{y_{1}}^{\pm 1}, {y_{2}}^{\pm 1}, \dots {y_{p^{n}}}^{\pm 1}]$ becomes an $R$-Hopf algebra by  
\begin{align*}
&\text{(a) the multiplication}: y_{s} \mapsto y_{s} \otimes y_{s} \  \text{for} \  1 \leq s \leq p^{n}; \ \ \ \ \ \ \ \ \ \ \ \ \ \ \ \ \ \ \ \ \ \ \ \ \ \ \ \ \ \ \ \ \ \ \ \ \ \ \ \ \\
&\text{(b) the unit}: y_{s} \mapsto 1  \  \text{for} \  1 \leq s \leq p^{n};\\
&\text{(c) the inverse}: y_{s} \mapsto \frac{1}{y_{s}} \  \text{for} \  1 \leq s \leq p^{n}.
\end{align*}
Then, $\xi$ is an $R$-Hopf algebras map. 

\vspace{3mm}

\begin{notation} In the coordinate ring $R[\Lambda, \frac{1}{\Lambda}][X_{1}, X_{T}, X_{T^{2}}, \dots X_{T^{p^{n}-1}}, {1}/{{D_{(\Lambda)}}}  ]$ of $U(\varGamma^{(\Lambda)}_{R[\Lambda, \frac{1}{\Lambda}]})$, we put
\[P_{(\Lambda)}( X^{2}_{T})=\frac{1}{\Lambda^{2}}\Bigr\{\frac{(X_{1}+\Lambda X_{T})^{2}}{X_{1}+2\Lambda X_{T}+ \Lambda^{2}X_{T^{2}}}-X_{1}    \Bigr\}.\]

For $3 \leq s \leq p^{n}-1$, we put

\[P_{(\Lambda)}( X^{s}_{T})=\frac{1}{\Lambda^{s}}\Bigr\{ \frac{(X_{1}+\Lambda X_{T})^{s}}{\displaystyle\sum_{k=0}^{s}{\binom{s}{k}}\Lambda^{k}X_{T^{k}}}-X^{s-1}_{1}-\sum_{k=2}^{s-1}{\binom{s}{k}\Lambda^{k}X_{1}^{s-k}  P_{(\Lambda)}( X^{k}_{T})  }   \Bigr\}.\]
\end{notation}

\vspace{3mm}

$U(\varGamma^{(\Lambda)}_{R[\Lambda, \frac{1}{\Lambda}]})$is a
cleft $\varGamma^{(\Lambda)}_{R[\Lambda, \frac{1}{\Lambda}]}$-torsor over $U(\varGamma^{(\Lambda)}_{R[\Lambda, \frac{1}{\Lambda}]})/\varGamma^{(\Lambda)}_{R[\Lambda, \frac{1}{\Lambda}]}$, and 
$U(\boldsymbol{\mu}_{p^{n}, R[\Lambda, \frac{1}{\Lambda}]})$ is also
a cleft $\boldsymbol{\mu}_{p^{n}, R[\Lambda, \frac{1}{\Lambda}]}$-torsor over $U(\boldsymbol{\mu}_{p^{n}, R[\Lambda, \frac{1}{\Lambda}]})/\boldsymbol{\mu}_{p^{n}, R[\Lambda, \frac{1}{\Lambda}]}$.
Therefore, by considering the isomorphism $U(\widetilde{\alpha^{(\Lambda)}})$, we obtain the following commutative diagram
\[
\xymatrix{
 {S(A_{\boldsymbol{\mu}, R[\Lambda, 1/\Lambda]})_{\Theta}}^{\mathrm{co}{A_{\boldsymbol{\mu}, R[\Lambda, 1/\Lambda]}}} \otimes A_{\boldsymbol{\mu}, R[\Lambda, 1/\Lambda]}\ar[r]^-{\Phi_{R[\Lambda, 1/\Lambda]}}\ar[d]_-{U(\widetilde{\alpha^{(\Lambda)}})_{1}^{\#} \otimes \widetilde{\alpha^{(\Lambda)}}^{\#}}\ar@{}[rd]|{}&S(A_{\boldsymbol{\mu}, R[\Lambda, 1/\Lambda]})_{\Theta} \ar[d]^-{U(\widetilde{\alpha^{(\Lambda)}})^{\#}}\\
 {S(A_{R[\Lambda, 1/\Lambda]}^{(\Lambda)})_{\Theta}}^{\mathrm{co}{A_{R[\Lambda, 1/\Lambda]}^{(\Lambda)}}} \otimes A_{R[\Lambda, 1/\Lambda]}^{(\Lambda)}\ar[r]_-{\Psi_{ R[\Lambda, 1/\Lambda]}}& S(A_{R[\Lambda, 1/\Lambda]}^{(\Lambda)})_{\Theta}.
}
\]

Here, $U(\widetilde{\alpha^{(\Lambda)}})_{1}^{\#}$ is the restriction of $U(\widetilde{\alpha^{(\Lambda)}})^{\#}$ of $S(A_{\boldsymbol{\mu}, R[\Lambda, 1/\Lambda]})_{\Theta}^{\mathrm{co}{A_{\boldsymbol{\mu}, R[\Lambda, 1/\Lambda]}}}$. 

For $2 \leq s \leq p^{n}-1$, there exist $a_{k} \in  S(A_{R[\Lambda, 1/\Lambda]}^{(\lambda)})_{\Theta}^{\mathrm{co}{A_{R[\Lambda, 1/\Lambda]}^{(\lambda)}}} \  (0 \leq k \leq p^{n}-1)$ such that
\[\Psi_{R[\Lambda, 1/\Lambda]}^{-1}({X}_{T}^{s})=\sum_{k=0}^{p^{n}-1}a_{k} \otimes T^{k}.\]

From the above commutative diagram, for $s=2$, 
\[a_{0}=  P_{(\Lambda)}( X^{2}_{T})=\frac{1}{\Lambda^{2}}\Bigr\{\frac{(X_{1}+\Lambda X_{T})^{2}}{X_{1}+2\Lambda X_{T}+ \Lambda^{2}X_{T^{2}}}-X_{1}    \Bigr\}.    \]
Therefore, by induction, for $3 \leq s \leq p^{n}-1$, 
\[a_{0}=P_{(\Lambda)}( X^{s}_{T})=\frac{1}{\Lambda^{s}}\Bigr\{ \frac{(X_{1}+\Lambda X_{T})^{s}}{\displaystyle\sum_{k=0}^{s}{\binom{s}{k}}\Lambda^{k}X_{T^{k}}}-X^{s-1}_{1}-\sum_{k=2}^{s-1}{\binom{s}{k}\Lambda^{k}X_{1}^{s-k}  P_{(\Lambda)}( X^{k}_{T})  }   \Bigr\}.\]

Furthermore, we obtain the following proposition. 
\vspace{3mm}

\begin{prop} For $2 \leq s \leq p^{n}-1$,

\[ P_{(\Lambda)}( X^{s}_{T}) \in R[\Lambda][X_{1}, X_{T}, X_{T^{2}}, \dots X_{T^{p^{n}-1}}, \frac{1}{{D_{(\Lambda)}}}].\]

\end{prop}
\vspace{3mm}
\begin{proof}
The cleaving map

\[\psi_{R[\Lambda, 1/\Lambda]} : A_{R[\Lambda, 1/\Lambda]}^{(\lambda)} \rightarrow   S(A_{R[\Lambda, 1/\Lambda]}^{(\lambda)})_{\Theta}\] 
is given by 
\[T^{s} \mapsto X_{T^{s}} \  \  0 \leq s \leq p^{n}-1.\]
Since
$\varDelta (T^{0})=T^{0} \otimes T^{0}$,  $\varepsilon(T_{0})=1$, 
we have
\[\psi_{R[\Lambda, 1/\Lambda]}(T^{0})\psi^{-1}_{R[\Lambda, 1/\Lambda]}(T^{0})=1.\]
Hence,
$\psi^{-1}_{R[\Lambda, 1/\Lambda]}(T^{0})=1/X_{1}$. 

Next, since $\varDelta(T)=T\otimes 1 + (1 + \Lambda T) \otimes T, \varepsilon(T)=0$,
we have
\[\psi_{R[\Lambda, 1/\Lambda]}(T)\psi^{-1}_{R[\Lambda, 1/\Lambda]}(1)+\psi_{R[\Lambda, 1/\Lambda]}(1+\Lambda T)\psi^{-1}_{R[\Lambda, 1/\Lambda]}(T)=0.\] 
Hence,
\[\psi^{-1}_{R[\Lambda, 1/\Lambda]}(T) \in R[\Lambda][X_{1}, X_{T}, X_{T^{2}}, \dots X_{T^{p^{n}-1}}, \frac{1}{{D_{(\Lambda)}}}].\]

If, for $0 \leq  k \leq s <p^{n}-1$, 

\[\psi^{-1}_{R[\Lambda, 1/\Lambda]}(T^{k}) \in R[\Lambda][X_{1}, X_{T}, X_{T^{2}}, \dots X_{T^{p^{n}-1}}, \frac{1}{{D_{(\Lambda)}}}]\]
then

\[\sum_{l=0}^{s+1}{\binom{s+1}{l}}\Big\{\sum_{k=0}^{l}{\binom{l}{k}}\Lambda^{k}X_{T^{k+s-l}}\Big\}\otimes \psi^{-1}_{R[\Lambda, 1/\Lambda]}(T^{l}) =0.\]

Therefore, we obtain

\[\psi^{-1}_{R[\Lambda, 1/\Lambda]}(T^{s+1}) \in R[\Lambda][X_{1}, X_{T}, X_{T^{2}}, \dots X_{T^{p^{n}-1}}, \frac{1}{{D_{(\Lambda)}}}].\]

Inductively, for $0 \leq m \leq p^{n}-1$, 
we obtain
\[\psi^{-1}_{R[\Lambda, 1/\Lambda]}(T^{m}) \in R[\Lambda][X_{1}, X_{T}, X_{T^{2}}, \dots X_{T^{p^{n}-1}}, \frac{1}{{D_{(\Lambda)}}}].\]

Let $\rho$ be the structure of right $A_{R[\Lambda, 1/\Lambda]}^{(\lambda)}$-comodule algebra, $S(A_{R[\Lambda, 1/\Lambda]}^{(\lambda)})$. For $0 \leq s \leq p^{n}-1$,
\[\rho(X^{s}_{T})=\sum_{k=0}^{s}{\binom{s}{k}}X^{k}_{T}(X_{1}+\Lambda X_{T})^{s-k} \otimes T^{s-k}\]
\[=\sum_{k=0}^{s}{\binom{s}{k}}X^{k}_{T}\Big\{\sum_{l=0}^{s-k}{\binom{s-k}{l}}X^{l}_{1}(\Lambda X_{T})^{s-k-l}\Big\} \otimes T^{s-k}.\]
Therefore,

\[P_{(\Lambda)}(X^{s}_{T})=\sum_{k=0}^{s}{\binom{s}{k}}X^{k}_{T}\Big\{\sum_{l=0}^{s-k}{\binom{s-k}{l}}X^{l}_{1}(\Lambda X_{T})^{s-k-l}\Big\} \otimes \psi^{-1}_{R[\Lambda, 1/\Lambda]}(T^{s-k})\] from the equation $(*)$ in the results of Doi and Takeuchi.

Hence, for $0 \leq s \leq p^{n}-1$,

\[ P_{(\Lambda)}( X^{s}_{T}) \in R[\Lambda][X_{1}, X_{T}, X_{T^{2}}, \dots X_{T^{p^{n}-1}}, \frac{1}{{D_{(\Lambda)}}}].\]
\end{proof}

\vspace{3mm}

Since $P_{(\Lambda)}( X^{s}_{T}) \in {S(A_{R[\Lambda, 1/\Lambda]}^{(\lambda)})_{\Theta}}^{\mathrm{co}{A_{R[\Lambda, 1/\Lambda]}^{(\lambda)}}}$ and $R[\Lambda, 1/\Lambda]$ is $R[\Lambda]$-flat, $P_{(\Lambda)}( X^{s}_{T}) \in {S(A_{R[\Lambda]}^{(\lambda)})_{\Theta}}^{\mathrm{co}{A_{R[\Lambda]}^{(\lambda)}}}$.

Moreover, from the equation $(*)$, there exist
\[a_{0}, a_{1}, \dots, a_{p^{n}-1} \in R[\Lambda][X_{1}^{\pm 1}, P_{(\Lambda)}( X^{2}_{T}), \dots,  P_{(\Lambda)}( X^{p^{n}-1}_{T}) ]  \] 
such that 
\[\Psi_{R[\Lambda, 1/\Lambda]}^{-1}({X}_{T}^{s})=\sum_{k=0}^{p^{n}-1}a_{k} \otimes T^{k}\]
for $2 \leq s \leq p^{n}-1$.

We have the following commutative diagram

\[
\xymatrix{
 {S(A_{R[\Lambda]}^{(\lambda)})_{\Theta}}^{\mathrm{co}{A_{R[\Lambda]}^{(\lambda)}}} \otimes A_{R[\Lambda]}^{(\lambda)} \ar[r]^-{\Psi_{R[\Lambda]}}\ar[d]_-{}\ar@{}[rd]|{}&S(A_{R[\Lambda]}^{(\lambda)})_{\Theta} \ar[d]^-{}\\
 {S(A_{R[\Lambda, 1/\Lambda]}^{(\lambda)})_{\Theta}}^{\mathrm{co}{A_{R[\Lambda, 1/\Lambda]}^{(\lambda)}}} \otimes A_{R[\Lambda, 1/\Lambda]}^{(\lambda)}\ar[r]_-{\ \ \Psi_{ R[\Lambda, 1/\Lambda]}}& S(A_{R[\Lambda, 1/\Lambda]}^{(\lambda)})_{\Theta}.
}
\]

Therefore, we obtain that $\Psi^{-1}_{R[\Lambda, 1/\Lambda]}(X^{s}_T)=\Psi^{-1}_{R[\Lambda]}(X^{s}_T)$ for $2 \leq s \leq p^{n}-1$.

Furthermore, by substituting $\lambda \in R$ for $\Lambda$, we obtain that

\[P_{(\lambda)}( X^{s}_{T}) \in {S(A_{R}^{(\lambda)})_{\Theta}}^{\mathrm{co}{A_{R}^{(\lambda)}}}\]

and there exist \[a_{0}, a_{1}, \dots, a_{ p^{n} -1} \in R[X_{1}^{\pm 1}, P_{(\lambda)}( X^{2}_{T}), \dots,  P_{(\lambda)} ( X^{p^{n}-1}_{T}) ]  \] 
such that
\[{\Psi}^{-1}(X_{T}^{s})=\sum_{k=0}^{p^{n}-1}a_{k} \otimes T^{k}\]

for $2 \leq s \leq p^{n}-1$.

\vspace{3mm}

\begin{prop}  $S(A_{R}^{(\lambda)})_{\Theta}^{\mathrm{co}A_{R}^{(\lambda)}}$ is the subalgebra of $S(A_{R}^{(\lambda)})_{\Theta}$ generated by the elements
\[\begin{split}
&X_{1}^{\pm 1}, X^{p^{n}}_{T}, P_{(\lambda)}(X^{2}_{T}), \dots, P_{(\lambda)}(X^{p^{n}-1}_{T}), \frac{1}{(X^{p^{n}}_{1}+ {\lambda}^{p^{n}}X^{p^{n}}_{T})},\\
& \frac{(\displaystyle\sum_{k=0}^{s} {\binom{s}{k}}\lambda^{k}X_{T^{k}})}{(X_{1}+\lambda X_{T})^{s}}, 2 \le s \le p^{n}-1. 
\end{split}\]
\end{prop}

\begin{proof} Let $B$ denote the subalgebra of $S(A_{R}^{(\lambda)})_{\Theta}$ generated by the elements
\[\begin{split}
&X_{1}^{\pm 1}, X^{p^{n}}_{T}, P_{(\lambda)}(X^{2}_{T}), \dots, P_{(\lambda)}(X^{p^{n}-1}_{T}), \frac{1}{(X^{p^{n}}_{1}+ {\lambda}^{p^{n}}X^{p^{n}}_{T})},\\
& \frac{(\displaystyle\sum_{k=0}^{s} {\binom{s}{k}}\lambda^{k}X_{T^{k}})}{(X_{1}+\lambda X_{T})^{s}}, 2 \le s \le p^{n}-1. 
\end{split}\]

First, we have \[B\subset S(A_{R}^{(\lambda)})_{\Theta}^{\mathrm{co}A_{R}^{(\lambda)}}.\]

Any nonnegative integer $m$ can be represented as $m=p^{n}k+r$ for some nonnegative integer $k$ and some nonnegative integer $r$ such that $0 \leq r <p^{n}$. 

Therefore, 
\[\Psi^{-1}(X^{m}_{T})=\Psi^{-1}(X^{p^{n}k+r}_{T})=\Big\{(X^{p^{n}}_{T})^{k} \otimes 1 \Big\} \Psi^{-1}(X^{r}_{T}).\]

Moreover, 

\[\Psi^{-1}\Big(\frac{1}{X_{1}+ \lambda X_{T}}\Big)=\Psi^{-1}\Big(\frac{({X_{1}+ \lambda X_{T}})^{p^{n}-1}}{({X_{1}+ \lambda X_{T}})^{p^n}}\Big)=\Big\{\frac{1}{({X_{1}+ \lambda X_{T}})^{p^n}} \otimes 1 \Big\}\Psi^{-1}\Big(({X_{1}+ \lambda X_{T}})^{p^{n}-1}\Big)\]
\[=\Bigr(\frac{1}{({X_{1}+ \lambda X_{T}})^{p^n}} \otimes 1 \Bigr)\Bigr\{\sum_{k=0}^{p^{n}-1}{\binom{p^{n}-1}{k}}\lambda^{k}(X^{p^{n}-k-1}_{1} \otimes 1)\Psi^{-1}(X^{k}_{T})\Bigr\}.\]

Hence, for any nonnegative integer $m$, 
there exist $a_{k} \in B$ for $0 \leq k \leq p^{n}-1$ such that
\[\Psi^{-1}\Big(\frac{1}{(X_{1}+\lambda X_{T})^{m}}\Big)=\sum_{k=0}^{p^{n}-1}(a_{k} \otimes 1)\Psi^{-1}(X^{k}_{T}).\]

For $2 \leq s \leq p^{n}-1$, 
\[\Psi^{-1}\Big(\frac{1}{\displaystyle\sum_{k=0}^{s}{\binom{s}{k}}\lambda^{k}X_{T^{k}}}\Big)=\Psi^{-1}\Big(\frac{(X_{1}+\lambda X_{T})^{s}}{\Big\{\displaystyle\sum_{k=0}^{s}{\binom{s}{k}}\lambda^{k}X_{T^{k}}\Big\}(X_{1}+\lambda X_{T})^{s}}\Big)\]
\[=\Big(X_{1}^{s-1}+\displaystyle\sum_{k=2}^{s-1}{\binom{s}{k}}\lambda^{k}X_{1}^{s-k}P_{\lambda}(X^{k}_{T})+ \lambda^{s}P_{(\lambda)}(X^{s}_{T}) \otimes 1\Big)\Psi^{-1}\Bigr(\frac{1}{(X_{1}+\lambda X_{T})^{s}}\Bigr).\]
Hence, 

for any nonnegative integer $m$, we find there exist $a_{k} \in B$ for $0 \leq k \leq p^{n}-1$ such that 
\[\Psi^{-1}\Big(\frac{1}{\Big\{\displaystyle\sum_{k=0}^{s}{\binom{s}{k}}\lambda^{k}X_{T^{k}}\Big\}^{m}}\Big)=\sum_{k=0}^{p^{n}-1}(a_{k} \otimes 1)\Psi^{-1}(X^{k}_{T}).\]
We have 
\[P_{(\lambda)}(X^{2}_{T})=\frac{-X_{1}X_{T^{2}}+Q_{2}}{X_{1}+2\lambda X_{T} + \lambda^{2}X_{T^{2}}},\]
where
\[Q_{2} \in R[X_{1}^{\pm 1}, X_{T}, \frac{1}{X_{1}+\lambda X_{T}}].\]
Since 
\[X_{1}+ \lambda^{2} P_{(\lambda)}(X^{2}_{T})=\frac{(X_{1}+\lambda X_{T})^{2}}{X_{1}+2 \lambda X_{T}+ \lambda^{2} X_{T^{2}}},\] we obtain that
\[X_{T^{2}}=\frac{Q_{2}-(X_{1}+2 \lambda X_{T}) P_{(\lambda)}(X^{2}_{T})}{X_{1}+ \lambda^{2} P_{(\lambda)}(X^{2}_{T})}.\]
Hence, for any nonnegative integer $m$, there exist $a_{k} \in B$ for $0 \leq k \leq p^{n}-1$
such that
\[\Psi^{-1}(X^{m}_{T^{2}})=\sum_{k=0}^{p^{n}-1}(a_{k} \otimes 1)\Psi^{-1}(X^{k}_{T}).\]

If there exist $a_{k} \in B$ for $0 \leq k \leq p^{n}-1$ such that
\[\Psi^{-1}(X_{T^{s}})=\sum_{k=0}^{p^{n}-1}(a_{k} \otimes 1)\Psi^{-1}(X^{k}_{T})\] 
 for $2 \leq s \leq p^{n}-2$, we have
\[P_{(\lambda)}(X^{s+1}_{T})=\frac{\Big\{-X^{s}_{1}-\displaystyle\sum_{k=2}^{s}{\binom{s+1}{k}\lambda^{k}X_{1}^{s+1-k}P_{(\lambda)}(X^{k}_{T})}\Big\}X_{T^{s+1}}+Q_{s+1}}{\displaystyle\sum_{k=0}^{s+1}{\binom{s+1}{k}}\lambda^{k}X_{T^{k}}},\]
where
\[Q_{s+1} \in R[X_{1}^{\pm 1}, X_{T}, \dots, X_{T^{s}}, \frac{1} {\displaystyle\prod_{l=1}^{s}\Big(\displaystyle\sum_{k=0}^{l}{\binom{l}{k}}\lambda^{k}X_{T^{k}}\Big)}]. \]
Hence, 
\[X^{s}_{1}+\sum_{k=2}^{s}{\binom{s+1}{k}\lambda^{k}X_{1}^{s+1-k}P_{(\lambda)}(X^{k}_{T})} + \lambda^{s+1} P_{(\lambda)}(X^{s+1}_{T})=\frac{(X_{1}+ \lambda X_{T})^{s+1}}{\displaystyle\sum_{k=0}^{s+1}{\binom{s+1}{k}}\lambda^{k}X_{T^{k}}},\]
we obtain that
\[X_{T^{s+1}}=\frac{-\Big\{\displaystyle\sum_{k=0}^{s}{\binom{s+1}{k}}\lambda^{k}X_{T^{k}}\Big\}P(X^{s+1}_{T})+Q_{s+1}}{X^{s}_{1}+\displaystyle\sum_{k=2}^{s}{\binom{s+1}{k}\lambda^{k}X_{1}^{s+1-k}P_{(\lambda)}(X^{k}_{T})} + \lambda^{s+1} P_{(\lambda)}(X^{s+1}_{T})}.\]
Therefore, there exist $a_{k} \in B$ for $0 \leq k \leq p^{n}-1$
such that
\[\Psi^{-1}(X^{m}_{T^{s+1}})=\sum_{k=0}^{p^{n}-1}(a_{k} \otimes 1)\Psi^{-1}(X^{k}_{T})\]
for any nonnegative integer $m$.

Hence, a homomorphism of $R$-modules
\[{\Psi}': B \otimes A_{R}^{(\lambda)} \rightarrow S(A_{R}^{(\lambda)})_{\Theta}  \]
defined by
\[b \otimes a \mapsto b \psi(a) \]
is bijective.

Considering the natural injection $s: B \rightarrow S(A_{R}^{(\lambda)})_{\Theta}^{\mathrm{co}A_{R}^{(\lambda)}}$,
we obtain the following commutative diagram. 

\[
  \xymatrix{
     B \otimes A_{R}^{(\lambda)}\ar[r]^{{ \ \ \ \  \Psi}'} \ar[d]_{s \otimes \text{Id}} &  S(A_{R}^{(\lambda)})_{\Theta}  \\
     {S(A_{R}^{(\lambda)})_{\Theta}}^{\mathrm{co}A_{R}^{(\lambda)}}  \otimes A_{R}^{(\lambda)}.  \ar[ur]^{\Psi} & 
  }
\]

Moreover, since $A_{R}^{(\lambda)}$ is faithfully flat over $R$, $s$ is bijective.

\end{proof}

\vspace{3mm}
\begin{prop} $S(A_{R}^{(\lambda)})_{\Theta}$ is isomorphic, as an $R$-algebra, to the polynomial algebra
$R[Z_{1}, Z_{T}, \cdots, Z_{T^{p^{n}-1}}]$ localized by the elements

\[\begin{split}
&Z_{1}, Z_{1}+\lambda Z_{T}, Z_{1}+ \lambda^{2}Z_{T^{{2}}}, \\
&Z_{1}^{s-1}+\sum_{k=2}^{s-1}\binom{s}{k}\lambda^{k}Z_{1}^{s-k}Z_{T^{k}}+\lambda^{s}Z_{T^{s}}, 3 \le s \le p^{n}-1.
\end{split}\]
\end{prop}
\begin{proof}

Let $A_{\mathcal{W}}$ denote the polynomial algebra
$R[Z_{1}, Z_{T}, \cdots, Z_{T^{p^{n}-1}}]$ localized by the elements
\[\begin{split}
&Z_{1}, Z_{1}+\lambda Z_{T}, Z_{1}+ \lambda^{2}Z_{T^{{2}}}, \\
&Z_{1}^{s-1}+\sum_{k=2}^{s-1}\binom{s}{k}\lambda^{k}Z_{1}^{s-k}Z_{T^{k}}+\lambda^{s}Z_{T^{s}}, 3 \le s \le p^{n}-1.
\end{split}\]

We define a homomorphism of $R$-algebras
\[\chi:S(A_{R}^{(\lambda)})_{\Theta} \rightarrow  A_{\mathcal{W}}\]
define by
\[\chi(X_{1})=Z_{1}, \chi(X_{T})=Z_{T},\]

\[\chi(X_{T^{2}})=\frac{\chi(Q_{2})-(Z_{1}+2 \lambda Z_{T}) Z_{T^{2}}}{Z_{1}+ \lambda^{2} Z_{T^{2}}},\]

\[\chi(X_{T^{s+1}})=\frac{-\Big\{\displaystyle\sum_{k=0}^{s}{\binom{s+1}{k}} \lambda^{k}\chi(X_{T^{k}})\Big\}Z_{T^{s+1}}+\chi(Q_{s+1})}{\chi(X^{s}_{1})+\displaystyle\sum_{k=2}^{s}{\binom{s+1}{k}\lambda^{k} (Z_{1}^{s+1-k}) Z_{T^{k}} + \lambda^{s+1}Z_{T^{s+1} }}} \]
for $2 \leq s \leq p^{n}-2$. This is well-defined.

Indeed, we obtain 
\[(Z_{1}+2\lambda Z_{T}+\lambda^{2}\chi(X_{T^{2}}))Z_{T^{2}}=\chi(Q_{2})-Z_{1}\chi(X_{T^{2}})\]
since \[(Z_{1}+\lambda^{2}Z_{T^{2}})\chi(X_{T^{2}})=\chi(Q_{2})-(Z_{1}+2 \lambda Z_{T})Z_{T^{2}}.\]

Moreover, by multiplying both sides by $\lambda^{2}$, 
\[(Z_{1}+2\lambda Z_{T}+\lambda^{2}\chi(X_{T^{2}}))\lambda^{2}Z_{T^{2}}=(Z_{1}+\lambda Z_{T})^{2}-Z_{1}(Z_{1}+2\lambda Z_{T}+\lambda^{2}\chi(X_{T^{2}})).\]
Therefore,
\[(Z_{1}+2\lambda Z_{T}+\lambda^{2}\chi(X_{T^{2}}))(Z_{1}+\lambda^{2}Z_{T^{2}})=(Z_{1}+\lambda Z_{T})^{2}.\]

Hence, \[\frac{1}{(Z_{1}+2\lambda Z_{T}+\lambda^{2}\chi(X_{T^{2}}))}=\frac{Z^{2}_{1}+\lambda^{2}Z_{T^{2}}}{(Z_{1}+\lambda Z_{T})^{2}}.\]

For $2 \leq s \leq p^{n}-2$, we obtain
\[\Big\{\sum_{k=0}^{s}{\binom{s+1}{k}} \lambda^{k}\chi(X_{T^{k}})+\lambda^{s+1}\chi(X_{T^{s+1}})\Big\}Z_{T^{s+1}}\]
\[=\chi(Q_{s+1})-\Big\{Z^{s}_{1}+\sum_{k=2}^{s}\binom{s+1}{k}\lambda^{k} (Z_{1}^{s+1-k}) Z_{T^{k}}\Big\}\chi(X_{T^{s+1}})\]
since \[\Big\{Z^{s}_{1}+\sum_{k=2}^{s}{\binom{s+1}{k}\lambda^{k} (Z_{1}^{s+1-k}) Z_{T^{k}} + \lambda^{s+1}Z_{T^{s+1}}}\Big\}Z_{T^{s+1}}\]\[=-\Big\{\sum_{k=0}^{s}{\binom{s+1}{k}} \lambda^{k}\chi(X_{T^{k}})\Big\}Z_{T^{s+1}}+\chi(Q_{s+1}).\]
Moreover, by multiplying both sides by $\lambda^{s+1}$, 
\[\Big\{\sum_{k=0}^{s}{\binom{s+1}{k}} \lambda^{k}\chi(X_{T^{k}})+\lambda^{s+1}\chi(X_{T^{s+1}})\Big\}\lambda^{s+1}Z_{T^{s+1}}\]
\[=(Z_{1}+\lambda Z_{T})^{s+1}-\Big\{Z^{s}_{1}+\sum_{k=2}^{s}\binom{s+1}{k}\lambda^{k} Z_{1}^{s+1-k} Z_{T^{k}}\Big\}\Big\{\sum_{k=0}^{s+1}\binom{s+1}{k}\lambda^{k}\chi(X_{T^{k}})\Big\}\]
Therefore,
\[\Big\{Z^{s}_{1}+\sum_{k=2}^{s}\binom{s+1}{k}\lambda^{k} (Z_{1}^{s+1-k}) Z_{T^{k}}+\lambda^{s+1}Z_{T^{s+1}}\Big\}\Big\{\sum_{k=0}^{s+1}\binom{s+1}{k}\lambda^{k}\chi(X_{T^{k}})\Big\}\]
\[=(Z_{1}+\lambda Z_{T})^{s+1}.\]
Hence,
\[\frac{1}{\displaystyle\sum_{k=0}^{s+1}\binom{s+1}{k}\lambda^{k}\chi(X_{T^{k}})} =\frac{(Z_{1}+\lambda Z_{T})^{s+1}}{Z^{s}_{1}+\displaystyle\sum_{k=2}^{s}\binom{s+1}{k}\lambda^{k} (Z_{1}^{s+1-k}) Z_{T^{k}}+\lambda^{s+1}Z_{T^{s+1}}}\]

Moreover, we define a homomorphism of $R$-algebras
\[\xi: A_{\mathcal{W}} \rightarrow S(A_{R}^{(\lambda)})_{\Theta} \]
defined by
\[Z_{1} \mapsto X_{1}, Z_{T} \mapsto X_{T},\]

\[Z_{T^{s}} \mapsto P_{(\lambda)}(X^{s}_{T})\]
for $2 \leq s \leq p^{n}-1$.

This is well-defined. Then $\xi \circ \chi = \mathrm{Id}$ and  $\chi \circ \xi = \mathrm{Id}$. Therefore $\chi$ is bijective.

\end{proof}

\vspace{3mm}
\begin{theorem} 
$S(A_{R}^{(\lambda)})_{\Theta}^{\mathrm{co}A_{R}^{(\lambda)}}$ is isomorphic as an $R$-algebra to the polynomial algebra
$R[Z_{1}, Z_{T}, \cdots, Z_{T^{p^{n}-1}}]$ localized by the elements

\[\begin{split}
&Z_{1}, Z_{1}^{p^{n}}+\lambda^{p^{n}} Z_{T}, Z_{1}+ \lambda^{2}Z_{T^{{2}}}, \\
&Z_{1}^{s-1}+\sum_{k=2}^{s-1}\binom{s}{k}\lambda^{k}Z_{1}^{s-k}Z_{T^{k}}+\lambda^{s}Z_{T^{s}}, 3 \le s \le p^{n}-1.
\end{split}\]
\end{theorem}

\begin{proof}
Let $A_{\mathcal{V}}$ denote the polynomial algebra
$R[Z_{1}, Z_{T}, \cdots, Z_{T^{p^{n}-1}}]$ localized by the elements
\[\begin{split}
&Z_{1}, Z_{1}^{p^{n}}+\lambda^{p^{n}} Z_{T}, Z_{1}+ \lambda^{2}Z_{T^{{2}}}, \\
&Z_{1}^{s-1}+\sum_{k=2}^{s-1}\binom{s}{k}\lambda^{k}Z_{1}^{s-k}Z_{T^{k}}+\lambda^{s}Z_{T^{s}}, 3 \le s \le p^{n}-1.
\end{split}\]

We define a homomorphism of $R$-algebras \[\omega: A_{\mathcal{V}} \rightarrow A_{\mathcal{W}} \]
defined by
\[Z_{1} \mapsto Z_{1}, Z_{T} \mapsto Z^{p^{n}}_{T},\]

\[Z_{T^{s}} \mapsto Z_{T^{s}},\]
for $2 \leq s \leq p^{n}-1$.
This is well-defined. Moreover, we have that $\omega$ is injective. Since $\mathrm{Im} (\xi  \circ \omega)= S(A_{R}^{(\lambda)})_{\Theta}^{\mathrm{co}A_{R}^{(\lambda)}}$, this proof is complete.

\end{proof}

\vspace{5mm}

\begin{theorem} We have the following commutative diagrams of group schemes over $R$ with exact rows

\[\begin{CD}
 0 @>>> {\varGamma}^{(\lambda)}_{R} @>e>>  {\mathcal G}^{(\lambda)}_{R}   @>{F^{n}}>> \mathcal{G}_{R}^{(\lambda^{p^{n}})} @>>> 0\\
 @.     @\vert     @VV\sigma_{1}V    @VV\tau_{1}V  \\
 0 @>>> {\varGamma}^{(\lambda)}_{R} @>i>>  U({\varGamma}^{(\lambda)}_{R})  @>q>>  U({\varGamma}^{(\lambda)}_{R}) / {\varGamma}^{(\lambda)}_{R} @>>> 0
\end{CD}\]
and
\[\begin{CD}
 0 @>>> {\varGamma}^{(\lambda)}_{R} @>i>>  U({\varGamma}^{(\lambda)}_{R})   @>q>> U({\varGamma}^{(\lambda)}_{R}) / {\varGamma}^{(\lambda)}_{R} @>>> 0\\
 @.     @\vert     @VV\sigma_{2}V    @VV\tau_{2}V  \\
 0 @>>> {\varGamma}^{(\lambda)}_{R} @>e>>  {\mathcal G}^{(\lambda)}_{R} @>{F^{n}}>> {\mathcal G}^{(\lambda^{p^{n}})}_{R} @>>> 0.
\end{CD}\]
\end{theorem}
\begin{proof}

\vspace{3mm}
The morphism $q: U(\varGamma^{(\lambda)}_{R}) \rightarrow U(\varGamma^{(\lambda)}_{R})/\varGamma^{(\lambda)}_{R}$ is given by the natural injection of $R$-Hopf algebras
$inj: S(A_{R}^{(\lambda)})_{\Theta}^{\mathrm{co}A_{R}^{(\lambda)}} \subset  S(A_{R}^{(\lambda)})_{\Theta}$.

Moreover, we define a homomorphism of $R$-group schemes $\sigma_{1}: \mathcal{G}_{R}^{(\lambda)} \rightarrow  U(\varGamma^{(\lambda)}_{R})$  
by
\[\sigma_{1}^{\#}: R[X_{1}, X_{T}, X_{T^{2}}, \dots, X_{T^{p^{n}-1}}, \frac{1}{D_{(\lambda)}}] \rightarrow R[T, \frac{1}{1+\lambda T} ]  \]
\[X_{T^{s}} \mapsto T^{s}  \ \text{for} \ 0 \leq s \leq p^{n}-1\]
This is well-defined.

Furthermore, $\Im(\sigma_{1}^{\#} \circ inj) \subset R[T^{p^{n}}, \frac{1}{1+ \lambda^{p^{n}}T^{p^{n}}}]$. Therefore, there exists a homomorphism of $R$-group schemes $\tau_{1}: \mathcal{G}_{R}^{(\lambda^{p^{n}})} \rightarrow  U(\varGamma^{(\lambda)}_{R})/\varGamma^{(\lambda)}_{R}$ such that the following diagram of group schemes 

\[\begin{CD}
 U(\varGamma^{(\lambda)}_{R})   @>{q}>> U(\varGamma^{(\lambda)}_{R})/\varGamma^{(\lambda)}_{R} \\
 @AA{\sigma_{1}}A  @AA{{\tau_{1}}}A \\
 \mathcal{G}_{R}^{(\lambda)}  @>>{F^{n}}> \mathcal{G}_{R}^{(\lambda^{p^{n}})} 
\end{CD}\]
is commutative.

Conversely, we define a homomorphism of $R$-group schemes $\sigma_{2}:  U(\varGamma^{(\lambda)}_{R}) \rightarrow \mathcal{G}_{R}^{(\lambda)}$ by
\[\sigma_{2}^{\#}: R[T, \frac{1}{1+\lambda T} ] \rightarrow  R[X_{1}, X_{T}, X_{T^{2}}, \dots, X_{T^{p^{n}-1}}, \frac{1}{D_{(\lambda)}}]\]
\[T \mapsto \frac{X_{T}}{X_{1}}.\]
This is well-defined.

Furthermore, $\Im(\sigma_{2}^{\#} \circ {F^{p^{n}}}^{\#}) \subset S(A_{R}^{(\lambda)})_{\Theta}^{\mathrm{co}A_{R}^{(\lambda)}}$. Therefore, there exists a homomorphism of $R$-group schemes $\tau_{2}:  U(\varGamma^{(\lambda)}_{R})/\varGamma^{(\lambda)}_{R} \rightarrow  \mathcal{G}_{R}^{(\lambda^{p^{n}})} $ such that the following diagram of group schemes 

\[\begin{CD}
 U(\varGamma^{(\lambda)}_{R})   @>{q}>> U(\varGamma^{(\lambda)}_{R})/\varGamma^{(\lambda)}_{R} \\
 @VV{\sigma_{2}}V  @VV{{\tau_{2}}}V \\
 \mathcal{G}_{R}^{(\lambda)}  @>>{F^{n}}> \mathcal{G}_{R}^{(\lambda^{p^{n}})} 
\end{CD}\]
is commutative.

Therefore, we have commutative diagrams with exact rows.

\end{proof}

\vspace{3mm}
\begin{corollary} Let $S$ be an $R$-scheme and $X$ a $\varGamma^{(\lambda)}_{R}$-torsor over $S$. Then, the $\varGamma^{(\lambda)}_{R}$-torsor $X$ is cleft if and only if the class $[X]$ belongs to $\Ker[H^{1}(S,\varGamma^{(\lambda)}_{R}) \rightarrow H^{1}(S, {\mathcal G}^{(\lambda)}_{R}) ]$.
\end{corollary}
\vspace{3mm}
\begin{proof} By Theorem 3.8, we obtain a commutative diagram of cohomology groups
\[\begin{CD}
 H^{1}(S,\varGamma^{(\lambda)}_{R}) @>{i}>> H^{1}(S,U(\varGamma^{(\lambda)}_{R})) \\
 @\vert  @VV{\sigma_{2}}V \\
 H^{1}(S,\varGamma^{(\lambda)}_{R}) @>>{e}> H^{1}(S,  {\mathcal G}^{(\lambda)}_{R})
\end{CD}\]
(cf. Demazure-Gabriel [2, Ch III, Prop.4.6]).
Hence, we obtain an implication\[\Ker[H^{1}(S,G^{(\lambda)}_{R})\rightarrow  H^{1}(S,U(\varGamma^{(\lambda)}_{R}))]\subset \Ker[H^{1}(S,\varGamma^{(\lambda)}_{R})\rightarrow H^{1}(S,  {\mathcal G}^{(\lambda)}_{R})].\]
Conversely, by Theorem 3.8, we obtain a commutative diagram of cohomology groups
\[\begin{CD}
 H^{1}(S,\varGamma^{(\lambda)}_{R}) @>{i}>> H^{1}(S,U(\varGamma^{(\lambda)}_{R})) \\
 @\vert  @AA{\sigma_{1}}A \\
 H^{1}(S,\varGamma^{(\lambda)}_{R}) @>>{e}> H^{1}(S,  {\mathcal G}^{(\lambda)}_{R}).
\end{CD}\]
Hence, we obtain the implication \[ \Ker[H^{1}(S,\varGamma^{(\lambda)}_{R})\rightarrow H^{1}(S,  {\mathcal G}^{(\lambda)}_{R})]\subset \Ker[H^{1}(S,\varGamma^{(\lambda)}_{R})\rightarrow  H^{1}(S,U(\varGamma^{(\lambda)}_{R}))]\].
\end{proof}

\vspace{3mm}

\begin{corollary} Let $R$ be an $\F_{p}$-algebra and $\lambda \in R$. Put $\varGamma^{(\lambda)}_{R}=\Ker[F^{n}:  {\mathcal G}^{(\lambda)}_{R} \rightarrow  {\mathcal G}^{(\lambda^{p^{n}})}_{R}]$ and $S=\Spec\,R$. Then, a $\varGamma^{(\lambda)}_{R}$-torsor $X$ over $S$ is cleft if and only if there exist morphisms 
$X \rightarrow  {\mathcal G}^{(\lambda)}_{R}$ and $S \rightarrow  {\mathcal G}^{(\lambda^{p^{n}})}_{R}$ such that the diagram 
\[\begin{CD}
 X @>>>  {\mathcal G}^{(\lambda)}_{R} \\
 @VVV  @VV{F^{n}}V \\
 S @>>>  {\mathcal G}^{(\lambda^{p^{n}})}_{R}
\end{CD}\]
is Cartesian.
\end{corollary}

\vspace{3mm}
\begin{corollary} Under the notation of Corollary 3.10, the following conditions are equivalent:

(a) Any $\varGamma^{(\lambda)}_{R}$-torsor over $R$ is cleft.

(b) The map $ {\mathcal G}^{(\lambda^{p^{n}})}_{R}(R) \rightarrow H^{1}(R, \varGamma^{(\lambda)}_{R})$ induced by the exact sequence
\[0 \longrightarrow \varGamma^{(\lambda)}_{R} \longrightarrow {\mathcal G}^{(\lambda)}_{R}  \longrightarrow {\mathcal G}^{(\lambda^{p^{n}})}_{R}  \longrightarrow 0\]
is surjective.

(c) The map $H^{1}(R, {\mathcal G}^{(\lambda)}_{R} )\rightarrow H^{1}(R, {\mathcal G}^{(\lambda^{p^{n}})}_{R} )$ induced by the Frobenius morphism $F^{n}: {\mathcal G}^{(\lambda)}_{R}  \rightarrow {\mathcal G}^{(\lambda^{p^{n}})}_{R} $
is injective.

\end{corollary}
\vspace{3mm}

\vspace{3mm}
We conclude this paper by presenting an example of non-cleft $\varGamma^{(\lambda)}_{R}$-torsors.

\vspace{1mm}
\begin{notation} Let $R$ be a commutative ring and $\lambda\in R$. Put $S=\Spec\,R$. As usual we denote by ${\mathcal G}^{(\lambda)}_{R}$ the affine group scheme defined with the coordinate ring $B=R[T,1/(1+\lambda T)]$ equipped with the comultiplication $T \mapsto T\otimes 1 + 1 \otimes T + \lambda T\otimes T$.

Now, let $a\in R$ and put $C=C_{(\lambda,a)}=R[X,1/(a+\lambda X)]$. A right coaction
\[\rho_C: C=R[X,\frac{1}{a+\lambda X}] \rightarrow C\otimes_R B=R[X,\frac{1}{a+\lambda X}]\otimes_R R[T,\frac{1}{1+\lambda T}]\]
is defined by
\[\rho_C(X)=X\otimes 1 +a\otimes T + \lambda X\otimes T.\]
We put $X_{(\lambda,a)}=\Spec\,C_{(\lambda,a)}$.

Furthermore, assume that $R$ is an $\F_p$-algebra. Then, the coordinate ring of $\varGamma^{(\lambda)}_{R}=\Ker[F:{\mathcal G}^{(\lambda)}_{R}\rightarrow {\mathcal G}^{(\lambda^p)}_{R}]$ is given by  $\Tilde{B}=R[T]/(T^p)$ with the comultiplication $T \mapsto T\otimes 1 + 1 \otimes T + \lambda T\otimes T$.

Now, let $c\in R$ and put $\Tilde{C}=\Tilde{C}_{(\lambda,a,c)}=R[X]/(X^p-c)$. A right coaction
\[\rho_{\Tilde{C}}: \Tilde{C}=R[X]/(X^p-c) \rightarrow \Tilde{C}\otimes_R \Tilde{B}=R[X]/(X^p-c)\otimes_R R[T]/(T^p)\]
is defined by
\[\rho_{\Tilde{C}}(X)=X\otimes 1 +a\otimes T + \lambda X\otimes T.\]
We put $\Tilde{X}_{(\lambda,a,c)}=\Spec\,\Tilde{C}_{(\lambda,a,c)}$.
\end{notation}
\vspace{3mm}

\begin{prop} Under the above notations, we have the following assertions.

\noindent(1) If $a$ is invertible in $R/(\lambda)$, then $X_{(\lambda,a)}$ is a ${\mathcal G}^{(\lambda)}_{R}$-torsor over $S$.

\noindent(2) Assume that $R$ is an $\F_p$-algebra. If $a^p+\lambda^pc$ is invertible in $R$, then $\Tilde{X}_{(\lambda,a,c)}$ is a $\varGamma^{(\lambda)}_{R}$-torsor over $S$. Moreover, the contracted product $\Tilde{X}_{(\lambda,a,c)} \vee^{\varGamma^{(\lambda)}_{R}} {\mathcal G}^{(\lambda)}_{R}$ is isomorphic to $X_{(\lambda,a)}$ as a right ${\mathcal G}^{(\lambda)}_{R}$-torsor.

\end{prop}
\vspace{1mm}
\begin{proof} (1) Put $C=C_{(\lambda,a)}=R[X,1/(a+\lambda X)]$. Then the $R$-algebra homomorphism $r: C\otimes_R C \rightarrow C \otimes_R B$ defined by
\[X\otimes 1\mapsto X\otimes 1,\ 1\otimes X\mapsto \lambda X\otimes T+X\otimes 1+a\otimes T\]
is bijective. The inverse of $r$ is given by
\[X\otimes 1\mapsto X\otimes1,\ 1\otimes T \rightarrow \frac{1\otimes X-X\otimes 1}{(a+\lambda X)\otimes 1}.\]

Hence, it remains to prove that $C$ is faithfully flat over $R$. First note that $C$ is flat over $R$ since $C$ is a fraction ring of the polynomial ring $R[X]$. Note that
\[C\otimes_R R/(\lambda)=R[X,\frac{1}{a+\lambda X}]\otimes_R R/(\lambda)=(R/(\lambda))[X]\]
since $a$ is invertible in $R/(\lambda)$. Conversely,
\[U\mapsto \frac{1-a}{\lambda}:R[\frac{1}{\lambda}][X,\frac{1}{a+\lambda X}]\rightarrow R[\frac{1}{\lambda}]\]
defines a section of the inclusion map $R[1/\lambda]\rightarrow R[1/\lambda][U, 1/(a+\lambda U)]$. These imply that $\Spec\,C\rightarrow \Spec\,R$ is surjective.

\vspace{1mm}
\noindent(2) Put $\Tilde{C}=\Tilde{C}_{(\lambda,a,c)}=R[X]/(X^p-c)$. We first remark that $a+\lambda X$ is invertible in $\Tilde{C}$ since $(a+\lambda X)^p=a^p+\lambda^p c$ is invertible in $R$. Therefore, the $R$-algebra homomorphism $r: C\otimes_R C \rightarrow C \otimes_R B$ defined by
\[X\otimes 1\mapsto X\otimes 1,\ 1\otimes X\mapsto \lambda X\otimes T+X\otimes 1+a\otimes T\]
is bijective. In fact, the inverse of $r$ is given by
\[X\otimes 1\mapsto X\otimes1,\ 1\otimes T \rightarrow \frac{1\otimes X-X\otimes 1}{(a+\lambda X)\otimes 1}.\]
Moreover, $B$ is a free $R$-module of rank $p$ and is faithfully flat over $R$.

Now, define an $R$-algebra homomorphism
\[\varphi:R[X]\longrightarrow \Tilde{C}\otimes_R B=R[X]/(X^p-c)\otimes_R R[T,\frac{1}{1+\lambda T}]\]
by
\[\varphi(X)=X\otimes 1 +a\otimes T + \lambda X\otimes T.\]
Noting that $\varphi(a+\lambda X)=(a+\lambda X)\otimes(1+\lambda T)$ and $(a+\lambda X)\otimes(1+\lambda T)\in (\Tilde{C}\otimes_R B)^{\times}$, we can conclude that $\varphi:R[X]\rightarrow \Tilde{C}\otimes_R B$ is extended uniquely to an $R$-algebra homomorphism
\[\varphi:C=R[X,\frac{1}{a+\lambda X}]\longrightarrow \Tilde{C}\otimes_R B=R[X]/(X^p-c)\otimes_R R[T,\frac{1}{1+\lambda T}].\]
It is readily seen that $\varphi:C\rightarrow \Tilde{C}\otimes_R B$ is a homomorphism of right $B$-comodules.

Furthermore, define a left action of $\varGamma^{(\lambda)}_{R}$ on $\Tilde{X}_{(\lambda,a,c)}\times_S {\mathcal G}^{(\lambda)}_{R}$ by $\gamma(x,g)=(x\gamma^{-1},\gamma g)$. In other words, the left action $\varGamma^{(\lambda)}_{R}\times_S (\Tilde{X}_{(\lambda,a,c)}\times_S {\mathcal G}^{(\lambda)}_{R})\rightarrow \Tilde{X}_{(\lambda,a,c)}\times_S {\mathcal G}^{(\lambda)}_{R}$ is defined by the left coaction
\begin{align*}
\Tilde{C}\otimes_R B\rightarrow &\Tilde{B}\otimes_R (\Tilde{C}\otimes_R B):\\
&X\otimes 1\mapsto 1\otimes X\otimes1+\frac{-T}{1+\lambda T}\otimes(a+\lambda X)\otimes1,\\
&1\otimes T\mapsto T\otimes 1\otimes1+1\otimes 1\otimes T+\lambda T\otimes 1\otimes T.
\end{align*}
Then, $\Im\,\varphi$ is contained in the invariant subalgebra $(\Tilde{C}\otimes_R B)^{\Tilde{B}}$. It follows that $\varphi$ induces an isomorphism  of ${\mathcal G}^{(\lambda)}_{R}$-torsors
\[\Tilde{X}_{(\lambda,a,c)} \vee^{\varGamma^{(\lambda)}_{R}} {\mathcal G}^{(\lambda)}_{R}=(\Tilde{X}_{(\lambda,a,c)}\times_S {\mathcal G}^{(\lambda)}_{R})/\varGamma^{(\lambda)}_{R} \overset{\sim}{\longrightarrow}X_{(\lambda,a)}.\]
\end{proof}
\vspace{3mm}
By combining Corollary 3.10 and Proposition 3.13, we obtain the following assertion.

\vspace{3mm}

\begin{corollary} Let $R$ be an $\F_p$-algebra and $\lambda,a,c\in R$. Assume that $a^p+\lambda^p c$ is invertible in $R$. Then, the following conditions are equivalent.

\noindent(a) The $\varGamma^{(\lambda)}_{R}$-torsor $\Tilde{X}_{(\lambda,a,c)}$ is cleft.

\noindent(b) The ${\mathcal G}^{(\lambda)}_{R}$-torsor $X_{(\lambda,a)}$ is trivial.

\noindent(c) There exists $b\in R$ such that $a+\lambda b$ is invertible in $R$.

\end{corollary}

\vspace{3mm}
Now, we return to the example presented in [15, Example 4.6]. Put
\[R=\F_p[X, Y, \frac{1}{X^p+Y^p+(X+1)^p Y}]\] 
and
\[\lambda=X+1,\ a=X+Y,\ c=Y.\]
Then, $a^p+\lambda^pc=X^p+Y^p+(X+1)^pY$ is invertible in $R$. Conversely, there does not exist $b(X,Y)\in R$ such that $(X+Y)+(X+1)b(X,Y)$ is invertible. It follows that the $\varGamma^{(\lambda)}_{R}$-torsor $\Tilde{X}_{(\lambda,a,c)}$ is not cleft.

\section{Acknowledgment}

This paper would not exist without the guidance of Professor Emeritus Noriyuki Suwa. The author sincerely thanks Professor Akira Masuoka for his valuable advice on the quotient group scheme. 
The author thanks Editage for the English language editing.

\end{document}